\documentclass[preprint,12pt,authoryear,5p]{elsarticle}

\usepackage{amssymb}
\usepackage{amsmath}

\usepackage[acronym]{glossaries}

\usepackage[inline]{enumitem}
\usepackage{ulem}
\usepackage[linesnumbered,ruled,vlined,commentsnumbered]{algorithm2e}

\SetCommentSty{mycommfont}

\usepackage{url}
\usepackage{booktabs}

\usepackage{tabularx}

\usepackage{todonotes}

\usepackage{cleveref}

\crefname{equation}{equation}{}
\crefrangelabelformat{equation}{(#3#1#4)\text{-}(#5#2#6)}

\newlength{\myleftlen}

\newcommand{\backup}{\hskip-\myleftlen\mkern-35mu}

\newacronym{ACM}{ACM}{Adsorption Cooling Machine}
\newacronym{OCP}{OCP}{Optimal Control Problem}
\newacronym{MINLP}{MINLP}{Mixed-Integer Non-Linear Program}
\newacronym{MILP}{MILP}{Mixed-Integer Linear Program}
\newacronym{HVAC}{HVAC}{Heating, Ventilation and Air Conditioning}
\newacronym{COP}{COP}{Coefficient Of Performance}
\newacronym{MPC}{MPC}{Model Predictive Control}
\newacronym{NMPC}{NMPC}{Non-Linear Model Predictive Control}
\newacronym{STCS}{STCS}{Solar Thermal Climate System}
\newacronym{HTS}{HTS}{High Temperature Storage}
\newacronym{HT}{HT}{High Temperature}
\newacronym{LTS}{LTS}{Low Temperature Storage}
\newacronym{LT}{LT}{Low Temperature}
\newacronym{MT}{MT}{Medium Temperature}
\newacronym{NLP}{NLP}{Non-Linear Program}
\newacronym{MIOCP}{MIOCP}{Mixed-Integer Optimal Control Problem}
\newacronym{MIOC}{MIOC}{Mixed-Integer Optimal Control}
\newacronym{CIA}{CIA}{Combinatorial Integral Approximation}
\newacronym{HD}{HD}{Heat Dissipation}
\newacronym{ODE}{ODE}{Ordinary Differential Equation}
\newacronym{DAE}{DAE}{Differential Algebraic Equation}
\newacronym{IPM}{IPM}{Interior Point Method}
\newacronym{SQP}{SQP}{Sequential Quadratic Programming}
\newacronym{LGPL}{LGPL}{Lesser General Public License}
\newacronym{MLD}{MLD}{Mixed Logical Dynamical}
\newacronym{MIQP}{MIQP}{Mixed-Integer Qua\-dra\-tic Program}
\newacronym{SUR}{SUR}{Sum-Up-Rounding}
\newacronym{LICQ}{LICQ}{Linear Independence Constraint Qualification}
\newacronym{IQR}{IQR}{Inter-Quartile Range}
\newacronym{BnB}{BnB}{Branch-and-Bound}
\newacronym{HP}{HP}{Heat Pump}
\newacronym{RT}{RT}{Recooling Tower}
\newacronym{CS}{CS}{Concrete Slab}
\newacronym{HR}{HR}{Heating Rod}
\newacronym{MHE}{MHE}{Moving Horizon Estimation}
\newacronym{EKF}{EKF}{Extended Kal\-man Filter}
\newacronym{RC}{RC}{Reinforcement Cage}
\newacronym{CC}{CC}{Conventional Control}
\newacronym{HM}{HM}{Heat Meter}

\makeglossaries

\glsunset{MPC} 

\journal{Applied Energy}

\begin{document}

\begin{frontmatter}

\title{A whole-year simulation study on nonlinear mixed-integer model predictive control for a thermal energy supply system with multi-use components}

\author[add1,add2]{Adrian B\"{u}rger\corref{cor1}}
\ead{adrian.buerger@hs-karlsruhe.de}
\cortext[cor1]{Corresponding author}
\author[add1,add3]{Markus Bohlayer}
\ead{markus.bohlayer@hs-karlsruhe.de}
\author[add1]{Sarah Hoffmann}
\author[add1]{Angelika Altmann-Dieses}
\ead{angelika.altmann-dieses@hs-karlsruhe.de}
\author[add1]{Marco Braun}
\ead{marco.braun@hs-karlsruhe.de}
\author[add2,add4]{Moritz Diehl}
\ead{moritz.diehl@imtek.uni-freiburg.de}

\address[add1]{Faculty of Management Science and Engineering, Karlsruhe University of Applied Sciences, Moltkestra\ss e~30, 76133~Karlsruhe, Germany.}
\address[add2]{Systems Control and Optimization Laboratory, Department of Microsystems Engineering (IMTEK), University of Freiburg, Georges-Koehler-Allee~102, 71990~Freiburg im Breisgau, Germany.}
\address[add3]{Friedrich-Alexander-Universit\"{a}t Erlangen-N\"{u}rnberg, School of Business and Economics, Lange Gasse 20, 90403 N\"{u}rnberg, Germany}
\address[add4]{Department of Mathematics, University of Freiburg, Ernst-Zermelo-Stra\ss e~1, 79104~Freiburg im Breisgau, Germany.}

\begin{abstract}

This work presents a whole-year simulation study on nonlinear mixed-integer Model Predictive Control (MPC) for a complex thermal energy supply system which consists of a heat pump, stratified water storages, free cooling facilities, and a large underground thermal storage. For solution of the arising Mixed-Integer Non-Linear Programs (MINLPs) we apply an existing general and optimal-control-suitable decomposition approach. To compensate deviation of forecast inputs from measured disturbances, we introduce a moving horizon estimation step within the MPC strategy. The MPC performance for this study, which consists of more than 50,000 real time suitable MINLP solutions, is compared to an elaborate conventional control strategy for the system. It is shown that MPC can significantly reduce the yearly energy consumption while providing a similar degree of constraint satisfaction, and autonomously identify previously unknown, beneficial operation modes.

\end{abstract}

\begin{keyword}
Model predictive control \sep energy systems \sep mixed-integer nonlinear programming

\end{keyword}

\end{frontmatter}

\section{Introduction}

The transition towards a sustainable energy system requires a substantial contribution in the buildings sector, which is responsible for approximately 40\,\% of the energy consumption in the EU \citep{Berardi2017}. Low-energy buildings enable efficient utilization of provided thermal energy and, in combination with thermal component activation, facilitate the deployment of free cooling and efficient use of heat pumps. An intelligent mix of heat sources and sinks in combination with smart storage management can increase the efficiency of these systems considerably, cf. \cite{Hesaraki2015}. However, such configurations might incorporate numerous operating modes and states, which makes it difficult to identify and implement conventional, rule-based control strategies that are able to determine the most efficient operating mode at any given time point.

During recent years, a strand of literature has built which shows that \gls{MPC} is well suited for control of energy systems. The incorporation of models and forecasts allows for efficient machinery utilization as well as situation-dependent, predictive control decisions, and renders \gls{MPC} well suited for systems governed by slow dynamics and/or large delays. Successful applications can be found for a variety of systems, such as heat pump systems \citep{Fischer2017,Baeten2017}, district heating systems \citep{Dainese2019} and microgrids \citep{Zachar2019, Iovine2019}.

Due to switching behavior of components such as heat pumps, combined heat and power plants etc., application of MPC for energy systems often requires the solution of mixed-integer optimization problems on real time suitable time scales. Therefore, according optimization problems are often formulated as either \glspl{MILP} as in \cite{Bianchini2019} and \cite{Lv2019}, or \glspl{MIQP} as in \cite{Vasallo2016} and \cite{Killian2018}, using linear(ized) models for system modeling.

However, the underlying processes are often characterized by nonlinear correlations, so that the utilization of nonlinear models for mixed-inte\-ger \gls{MPC} can improve the accuracy of system descriptions, and with this, the quality of control decisions. The solution of the arising \glspl{MINLP} on suitable time scales remains a challenging task though, which for non-trivial problems can only be achieved approximately, cf. \cite{Schweiger2017}, \cite{Dias2018}, \cite{Kuboth2019}.

Within this work, we show a successful application of a general, optimal-control-suitable approach for approximate \gls{MINLP} solution by \cite{Sager2009}, \cite{Sager2011a} within a whole-year simulation study for mixed-integer \gls{MPC} of the nonlinear model of a complex thermal energy supply system. The system consists of a \gls{HP}, stratified water storages, free cooling facilities, and a large underground thermal storage. Model nonlinearities arise, amongst others, from consideration of mass flow rates as continuous controls of the system as well as nonlinear \gls{COP} computations for the \gls{HP}, which further introduces switching behavior due to minimal part load. Since the quality of utilized forecasts highly influences performance (cf. \cite{Oldewurtel2012}), we introduce a \gls{MHE} step for online correction of forecasts within the \gls{MPC} strategy to compensate deviation of forecasted inputs from measured disturbances. The performance of the \gls{MPC} strategy applied for this study, which consists of more than 50,000 real time suitable \glspl{MINLP} solutions, is compared to an elaborate, conventional control strategy for the system.

The contributions of this paper are as follows. First, the use of techniques for reduction of nonlinearities and elimination of discontinuities on modeling of real-life sized energy systems for later application within derivative-based optimization methods are demonstrated. Then, the developed models are applied within an integrated approach for simulation of a complete \gls{MPC}-\gls{MHE} loop for a nonlinear switched thermal energy system. The methods applied for solving the arising \glspl{MINLP} are general and suitable not only for simulation, but also for real-time control applications of physical systems. Finally, it its shown within a case study that the \gls{MPC} is capable not only to significantly reduce the yearly energy consumption while providing a similar degree of constraint satisfaction, but also to identify beneficial operation modes which have not been identified before by the designers of the conventional controller.

The remainder of this paper is organized as follows. The system subject to this study is introduced in Section~\ref{sec:system_description}, while Section~\ref{sec:modeling_approach} presents the chosen modeling approach. The elaborate conventional control strategy is introduced in Section~\ref{sec:conventional_control} and the \gls{MPC} strategy in Section~\ref{sec:mpc_strategy}. The performance of the applied control strategies is compared within a simulation study in Section~\ref{sec:controller_performance_comparison}, followed by a conclusion and outlook in Section~\ref{sec:conclusion}.

\section{System description}
\label{sec:system_description}

In the following, a description of the system subject to this study is given by introducing its components, function principles, operation modes and limits.

\begin{figure*}[tbp]
\centering
\includegraphics[scale=0.84]{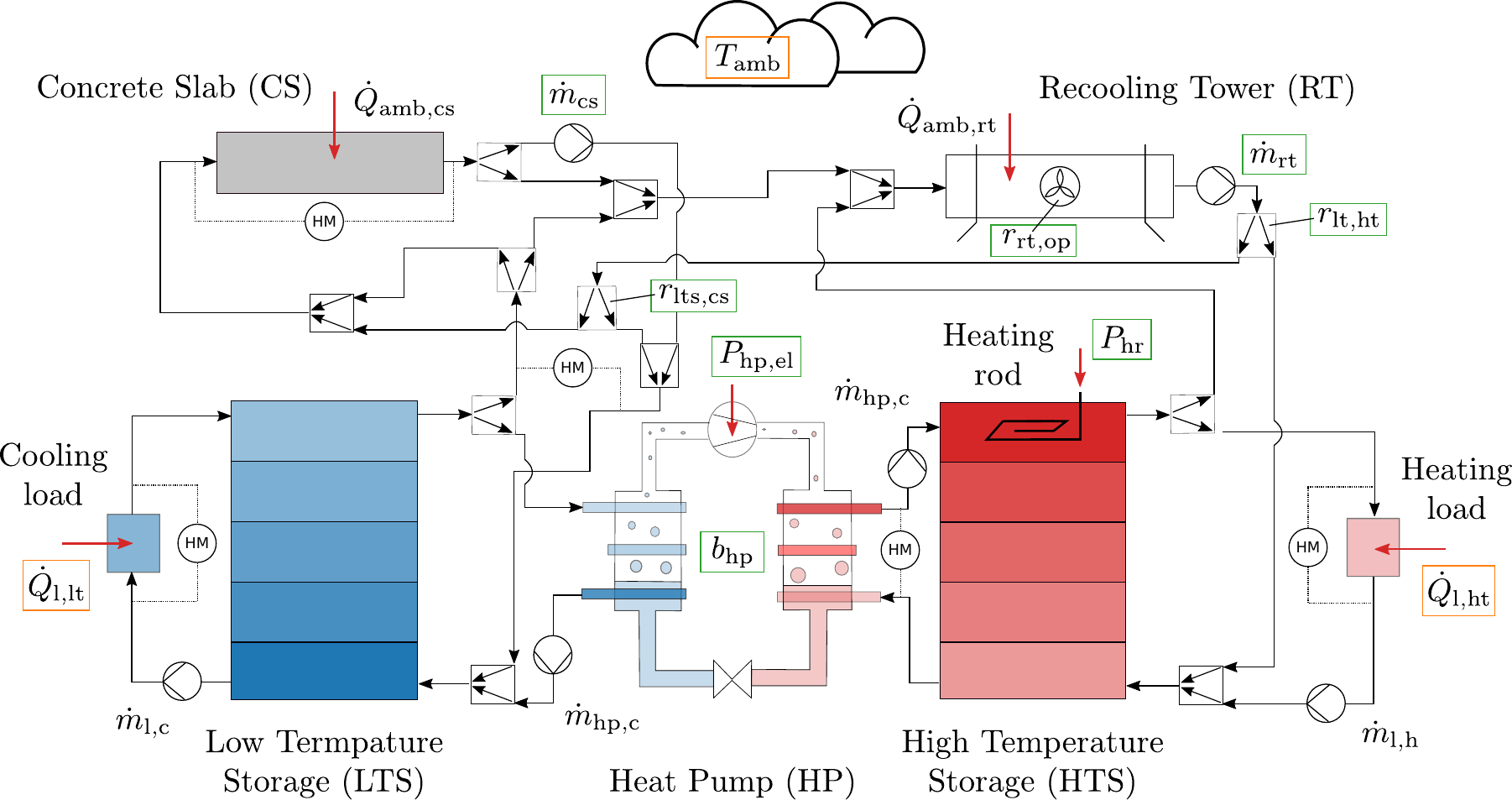}
\caption{Schematic depiction of the thermal energy supply system considered within this study: green boxes indicates controllable quantities, orange boxes indicate time-varying parameters determined by external factors. Red arrows indicate the main energy flows across the system boundaries.}
\label{fig:system_components}
\end{figure*}

\subsection{System components}

This study is inspired by an innovative thermal energy supply system for an office building under construction in Karls\-ruhe, Germany. The system variant and model used within this work originated from a research project of \textit{Ottensmeier Ingenieure GmbH, Paderborn} and \textit{Karls\-ruhe University of Applied Sciences} focused on the dimensioning of system components and control design for the physical system.

\begin{figure}[tbp]
\centering
\includegraphics[width=\columnwidth]{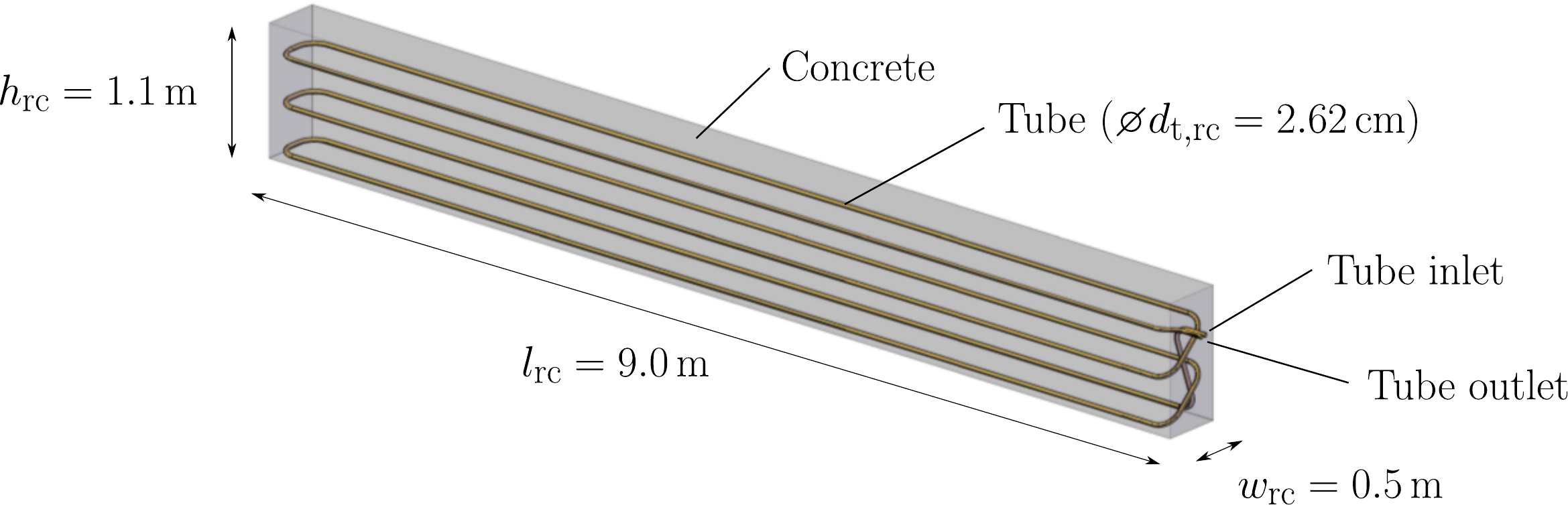}
\caption{Schematic depiction of one of the reinforcement cages that build the underground concrete slab.}
\label{fig:reinforcement_cage}
\end{figure}

An over\-view of the system considered within this study and the possible interactions between its components is given in Figure~\ref{fig:system_components}. One key component of the system is a \gls{CS} which is located below surface level of the building. This component, which assembles the underground parking deck of the building, is thermally activated through encased tubes within the concrete. These tubes with diameter $d_\mathrm{t,rc} = 2.62\,\mathrm{cm}$ are embedded as $n_\mathrm{cs} = 306 $ identical \glspl{RC}, exemplary shown in Figure~\ref{fig:reinforcement_cage}. Each cage is $w_\mathrm{rc} = 0.5\,\mathrm{m}$ wide, $h_\mathrm{rc} = 1.1\,\mathrm{m}$ high, $l_\mathrm{rc} = 9.0\,\mathrm{m}$ long and contains three layers of pipe loop with a total tube length of 58\,m per cage. Hydraulically, all cages are connected in parallel. The bottom of the \gls{CS} rests on soil, the top of the \gls{CS} is exposed to the air in the underground car park.

Further, two stratified water storages are installed: a \gls{LTS} with a volume of $V_\mathrm{lts} = 10.0\,\mathrm{m}^3$ and a \gls{HTS} with a volume of $V_\mathrm{hts} = 5.0\,\mathrm{m}^3$. Connected to these is a \gls{HP} with a maximum electric compressor power of $P_\mathrm{hp,el,max} = 60\,\mathrm{kW}$. Additionally, the \gls{HTS} is equipped with electrical \glspl{HR} with total maximum heating power $P_\mathrm{hr,max} = 400\,\mathrm{kW}$ intended as backup heating devices.\footnote{Within this study, \glspl{HR} are introduced as backup heating devices for simplicity. Alternatives could be, e.\,g., a connection to a district heating system.} Outside of the building, a \gls{RT} with total heat exchange area $A_\mathrm{rt} = 2145.0\,\mathrm{m}^2$, fluid volume $V_\mathrm{rt} = 0.208\,\mathrm{m}^3$ and heat transfer coefficient $\alpha_\mathrm {rt} = 21.6\,\mathrm{W}/\mathrm{K}$ is installed and connected to the indoor components.\footnote{For simplicity, usual separations of frost-proof outdoor circuits from indoor circuits is neglected within this study.} While the heating and cooling demands of the building on the one hand strongly depend on the current ambient temperature and solar irradiation, additional constant cooling demands arise from a number of server computers operated in the building. 

Multiple temperature sensors are installed inside \gls{CS}, storages and \gls{RT} as well as at their corresponding inlets and outlets. Additionally, five \glspl{HM} are installed in selected positions to enable monitoring and improved control decisions.

\subsection{Interaction of components}

The \gls{LTS} is the heat source for the \gls{HP}. When the \gls{HP} is active, water is extracted from the top of the \gls{LTS} and returned to the bottom of the storage with reduced temperature. From the bottom of the \gls{LTS}, water is extracted to cover cooling loads, resulting in water with increased temperature returning to the top of the storage.

Accordingly, the \gls{HTS} is the heat sink of the \gls{HP}, which extracts water from the bottom and returns it to the top of the storage with increased temperature. From the top of the \gls{HTS}, water is extracted to cover heating loads, resulting in water with reduced temperature returning to the bottom of the storage. The temperature of the \gls{HTS} can additionally be increased using the \glspl{HR}.
 
Medium from both storages can be supported to the outdoor \gls{RT} for heat exchange with the ambient air. Depending on current ambient temperature and operation level of the \gls{RT}, this allows to either cool down the \gls{LTS} in free cooling mode, or to heat up the \gls{LTS} when acting as a heat source for the \gls{HP}. On the other hand, the \gls{RT} can be utilized for cooling of the \gls{HTS} which enables the \gls{HTS} to act as recooling for the \gls{HP}.

Further, water from the \gls{LTS} can be transferred through the \gls{CS} to either reduce or increase the temperature inside of the \gls{LTS}, so that the \gls{CS} acts as both heat source and heat sink for the \gls{LTS}. Additionally, a splitter installed behind the \gls{RT} can be utilized to affect how much of the medium flowing through the \gls{RT} is led to the \gls{LTS} and how much is led into the \gls{CS}.

\subsection{Controls, operation limits and parameters}

Within this study, we assume that the mass flow $\dot{m}_\mathrm{rt}$ of the \gls{RT} pump and the mass flow $\dot{m}_\mathrm{cs}$ of the \gls{CS} pump can be chosen directly and on a continuous scale. Control $r_\mathrm{rt,op}$ sets the operation level of the \gls{RT}. The control $r_\mathrm{lt,ht}$ depicts a mixing valve which influences how much of the mass flow through the \gls{RT} is supported towards the \gls{LT} side of the system, while the remaining fraction is supported towards the \gls{HT} side. Control $r_\mathrm{lts,cs}$ is an additional mixing valve to determine how much of the medium flowing from the \gls{RT} to the \gls{LT} side is supported to the \gls{LTS}, while the remaining fraction is supported to the \gls{CS}.

The \gls{HP} can be activated using a binary switch $b_\mathrm{hp}$. Once the machine is turned on, it is assumed to operate with minimum part load $P_\mathrm{hp,el,min} = 2.0\,\mathrm{kW}$. From that, the electrical power $P_\mathrm{hp,el}$ of the \gls{HP} can be modulated continuously. The power of the \glspl{HR} can be chosen on a continuous scale. A summary of all controls is given in Table~\ref{tbl:controls}.

\begin{table}[tb]
\renewcommand*{\arraystretch}{1.3}
\centering
\scalebox{0.8}{%
\begin{tabular}{llll}
\hline
{Control} & {Type} & {Minimum} & {Maximum} \\
\hline
{$\dot{m}_\mathrm{rt}$} & {continuous} & {$0.0\,\mathrm{kg/s}$} & {$200.0\,\mathrm{kg/s}$} \\
{$\dot{m}_\mathrm{cs}$} & {continuous} & {$0.0\,\mathrm{kg/s}$} & {$75.0\,\mathrm{kg/s}$} \\
{$r_\mathrm{rt,op}$} & {continuous} & {0.0} & {1.0} \\
{$r_\mathrm{lt,ht}$} & {continuous} & {0.0} & {1.0} \\
{$r_\mathrm{lts,cs}$} & {continuous} & {0.0} & {1.0} \\
{$b_\mathrm{hp}$} & {binary} & {0} & {1} \\
{$P_\mathrm{hp,el}$} & {continuous} & {$2\,\mathrm{kW}$} & {$60\,\mathrm{kW}$} \\
{$\dot{Q}_\mathrm{hr}$} & {continuous} & {$0\,\mathrm{kW}$} & {$400\,\mathrm{kW}$} \\
\hline
\end{tabular}%
}
\caption{System controls and control boundaries.}
\label{tbl:controls}
\end{table}

The bottom temperature of the \gls{LTS} must not go below $3.0\,^\circ\mathrm{C}$ or exceed $18.0\,^\circ\mathrm{C}$, the top temperature of the \gls{HTS} must not go below $32.0\,^\circ\mathrm{C}$ or exceed $40.0\,^\circ\mathrm{C}$. The lower boundary for \gls{CS} temperatures is defined as $2.5\,^\circ\mathrm{C}$, the upper boundary as $40.0\,^\circ\mathrm{C}$. The maximum mass flow through the \gls{CS} is limited to a total of $\dot{m}_\mathrm{cs,max} = 75\,\mathrm{kg/s}$, the maximum mass flow through the \gls{RT} to a total of $\dot{m}_\mathrm{rt,max} = 200\,\mathrm{kg/s}$.

The ambient temperature $T_\mathrm{amb}$, heating demands $\dot{Q}_\mathrm{l,ht}$ and cooling demands $\dot{Q}_\mathrm{l,lt}$ enter the model as external time varying parameters.

\section{Modeling approach}
\label{sec:modeling_approach}

The models used within this study are nonlinear models based on mass and energy balances. If not stated explicitly, energy losses of components to their surroundings are neglected. Material and media is assumed incompressible and with constant specific heat capacities $c$ and densities $\rho$.

In the following, two model versions of the system will be used. For \gls{MPC} and \gls{MHE} computations, a model with reduced complexity adapted to fulfill all necessary requirements regarding differentiability for use with\-in derivative-based optimization methods, cf.~\cite{Biegler2010}, is utilized, further referred to as \textit{\gls{MPC} model} $f_\mathrm{mpc}$. For simulation of the conventional controller as well as for simulation of the \gls{MPC}-generated controls applied to the system, a more detailed model is used, further referred to as \textit{simulation model} $f_\mathrm{sim}$. Both models are systems of \glspl{ODE}, while the simulation model contains a total of $n_{x,\mathrm{sim}} = 199$ states and the \gls{MPC} model contains $n_{x,\mathrm{mpc}} = 36$ states. The models $f_\mathrm{sim}$ and $f_\mathrm{mpc}$ are supplied with different datasets for the external time-varying parameters, where the set $\{\dot{Q}_\mathrm{l,ht,fc},\,\dot{Q}_\mathrm{l,lt,fc},\,T_\mathrm{amb,fc}\}$ corresponds to forecast data used for \gls{MPC} computations and the set $\{\dot{Q}_\mathrm{l,ht,m},\,\dot{Q}_\mathrm{l,lt,m},\,T_\mathrm{amb,m}\}$ corresponds to measured data acting on the simulation model.

In the upcoming sections, the common modeling aspects for the several system components are given as well as a description of the differences between \gls{MPC} model and simulation model.

\subsection{Concrete slab}

\begin{figure}[tbp]
\centering
\includegraphics[width=0.7\columnwidth]{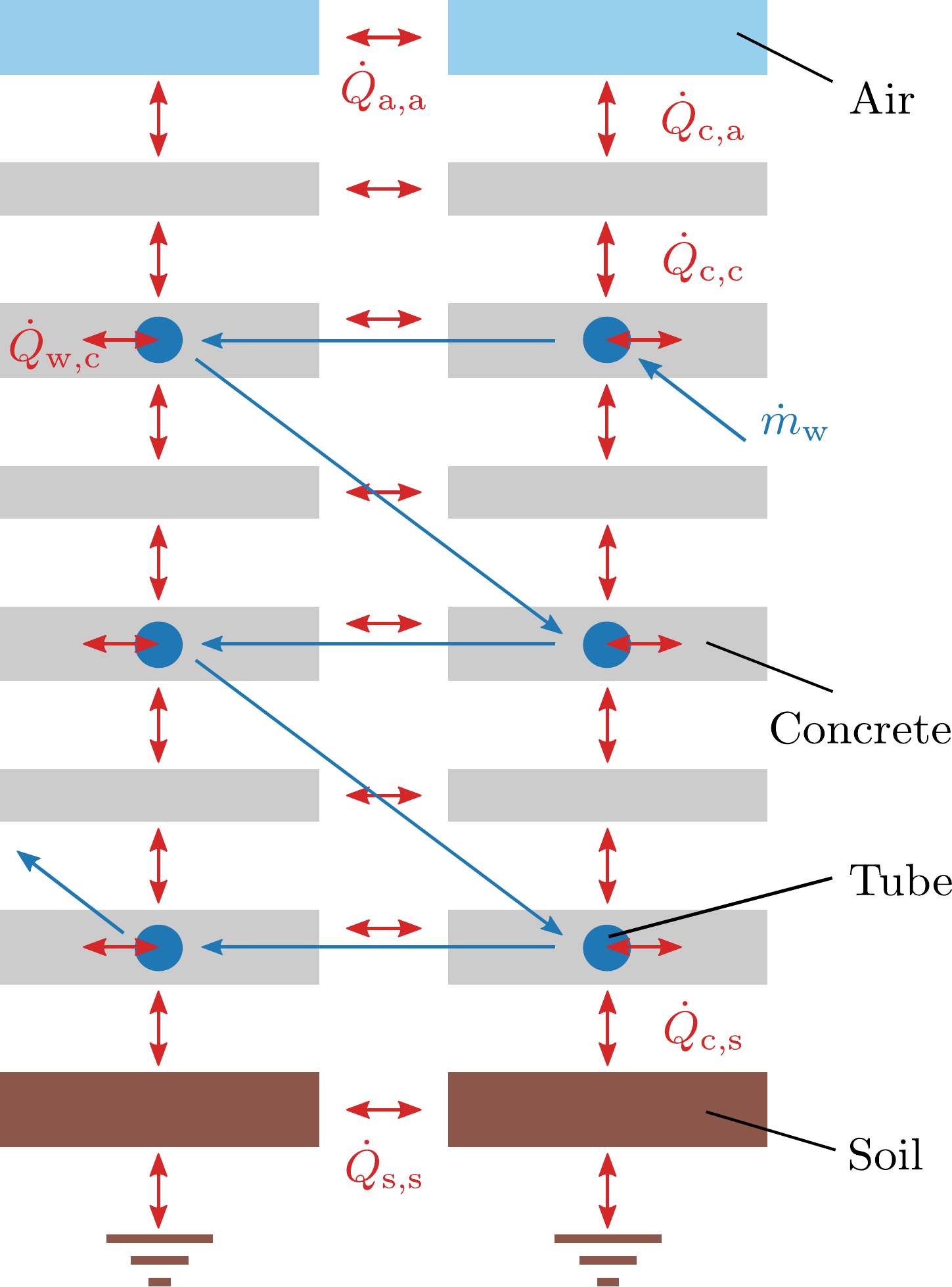}
\caption{Schematic depiction of the several balance volumes of the reinforcement cage model.}
\label{fig:cage_modeling}
\end{figure}

Since all \glspl{RC} within the \gls{CS} are connected in parallel, we can assume the temperature development within each \gls{RC} identical. Therefore, the \gls{CS} can be represented by a single cage that is supplied by a corresponding fraction of the total mass flow through the \gls{CS} and scaling the power according to the total number of cages. For modeling the temperature distribution inside the cage, its volume is discretized into a number of balance volumes of fluid and concrete. Additional volumes of soil and air build the cage's surrounding.  A schematic depiction showing the modeling approach on the discretization along the height and width of the \gls{RC} is given in Figure~\ref{fig:cage_modeling}.

The tube inside the cage is modeled by a number of fluid-filled volumes that exchange media and heat with neighboring volumes. The temperature $T_\mathrm{cs,w}$ within such a volume is calculated as
\begingroup
\allowdisplaybreaks
\begin{align}
\begin{split}
\dot{T}_\mathrm{cs,w}(t) =  {} & (\rho_\mathrm{w} V_\mathrm{cs,w})^{-1} \Bigl( \dot{m}_\mathrm{cs,t}(t)  (T_\mathrm{cs,w,pre}(t)  \\
 &  - T_\mathrm{cs,w}(t)) - c_\mathrm{w}^{-1}\dot{Q}_\mathrm{cs,w,c}(t) \Bigr)
 \end{split} \\ 
 \begin{split}
\dot{m}_\mathrm{cs,t}(t) =  {} & n_\mathrm{cs}^{-1} \Big(\dot{m}_\mathrm{cs}(t)  \\
& + r_\mathrm{lt,ht} (1.0 - r_\mathrm{lts,cs}(t)) \dot{m}_\mathrm{rt}(t) \Big)
\end{split} \label{eq:mdot_tube_sim_model} \\
\dot{Q}_\mathrm{cs,w,c}(t) = {} & \alpha_\mathrm{cs,w,c} A_\mathrm{cs,w,c} (T_\mathrm{cs,w}(t) - T_\mathrm{cs,c}(t))
\end{align}
\endgroup
with $V_\mathrm{cs,w}$ the constant fluid volume and $A_\mathrm{cs,w,c}$ the heat exchange area between fluid and concrete (i.\,e., the tube surface within the volume) depending on the discretization, $\alpha_\mathrm{cs,w,c}$ a constant heat transfer coefficient and $T_\mathrm{cs,c}$ the temperature of the surrounding concrete. For the volume at the inlet of the \gls{CS}, $T_\mathrm{cs,w,pre}$ results from the mixing of the \gls{RT} outlet temperature $T_\mathrm{rt}$ and the \gls{LTS} outlet temperature $T_\mathrm{lts}$ according to their mass flow ratios, for the upcoming volumes it is the temperature of the preceding volume.

Concrete within the cage is modeled by a number of solid volumes that exchange only heat with neighboring volumes. The temperature within a concrete volume $T_\mathrm{cs,c}$ is calculated as 
\begin{equation}
\begin{aligned}
\dot{T}_\mathrm{cs,c}(t) = {} & (\rho_\mathrm{c} V_\mathrm{cs,c} c_\mathrm{c})^{-1} \left( \dot{Q}_\mathrm{cs,c,s}(t) + \dot{Q}_\mathrm{cs,c,a}(t) \right. \\
& \left. + \dot{Q}_\mathrm{cs,w,c}(t) + \textstyle \sum_{i} \dot{Q}_{\mathrm{cs,c,c},i}(t) \right)
\end{aligned}
\end{equation}
where $\dot{Q}_\mathrm{cs,c,s}$ is the heat flowing from the center of the volume towards a neighboring soil volume (if existing), $\dot{Q}_\mathrm{cs,c,a}$ towards a neighboring air volume (if existing) and $\dot{Q}_{\mathrm{cs,c,c},i}$ towards the $i$th neighboring concrete volume. Accordingly, the temperature $T_\mathrm{cs,s}$ within a soil volume is calculated as
\begin{equation}
\dot{T}_\mathrm{cs,s}(t) = \frac{\dot{Q}_\mathrm{cs,s,c}(t) + \textstyle \sum_{i} \dot{Q}_{\mathrm{cs,s,s},i}(t)}{\rho_\mathrm{s} V_\mathrm{cs,s}c_\mathrm{s}}
\end{equation}
and the temperature $T_\mathrm{cs,a}$ within an air volume as
\begin{equation}
\dot{T}_\mathrm{cs,a}(t) = \frac{\dot{Q}_\mathrm{cs,a,c}(t) + \textstyle  \sum_{i} \dot{Q}_{\mathrm{a,a},i}(t)}{\rho_\mathrm{a} V_\mathrm{cs,a} c_\mathrm{a}}.
\end{equation}

Heat exchange between solid volumes is modeled depending on the adjoining heat exchange areas of volumes, the distance from the volume center to these areas and the thermal conductivity of the contained material, so that the current heat flow $\dot{Q}_\mathrm{x,y}$ between two solid volumes with temperatures $T_\mathrm{x}$ and $T_\mathrm{y}$ is given by
\begin{equation}
\dot{Q}_\mathrm{x,y}(t) = A \lambda_x \lambda_y \frac{T_\mathrm{x}(t) - T_\mathrm{y}(t)}{\Delta_x \lambda_y  + \Delta_y \lambda_x},
\end{equation}
with $\lambda_\mathrm{\{x,y\}}$ the thermal conductivities, $A$ the area of the adjoining surfaces and $\Delta_\mathrm{\{x,y\}}$ the distances of the center from the adjoining surface of a volume.

While for the MPC model the discretization is deducted only along the height and width of the cage as shown shown in Figure~\ref{fig:cage_modeling}, additional discretization using $n_l=5$ volumes along the length of the cage has been applied for the simulation model. Furthermore, within the \gls{MPC} model, the fraction of the mass flow from the \gls{RT} through the \gls{CS} is introduced as an individual control variable $\dot{m}_\mathrm{rt,cs}$ such that \eqref{eq:mdot_tube_sim_model} is replaced by
\begin{equation}
\dot{m}_\mathrm{cs,t}(t) = n_\mathrm{cs}^{-1} (\dot{m}_\mathrm{cs}(t) + \dot{m}_\mathrm{rt,cs}(t))
\end{equation}
to reduce nonlinearity of the \gls{MPC} model.

\subsection{Storages}

Both the \gls{LTS} and \gls{HTS} are modeled by a number of water-filled balance volumes. In both $f_\mathrm{sim}$ and $f_\mathrm{mpc}$, a discretization along the height of the storages using $n_\mathrm{lts} = n_\mathrm{hts} = 5$ volumes is deducted to depict temperature stratifications inside each storage, while the mass exchange and corresponding flow directions between the several volumes depend on the current mass flows at the storages in- and outlets, cf. \citet{Streckiene2011}.

For the simulation model, the energy balance that determines the temperature $T_{\mathrm{lts},1}$ of the \gls{LTS} bottom volume is formulated as
\begin{equation}
\begin{aligned}
\dot{T}_{\mathrm{lts},1}(t) =  {} & (\rho_\mathrm{w} V_\mathrm{lts})^{-1} \Bigl( \dot{m}_\mathrm{cs}(t) T_{\mathrm{cs,w},-1}(t)  \\
 & + r_\mathrm{lt,ht}(t) r_\mathrm{lts,cs}(t) \dot{m}_\mathrm{rt}(t) T_{\mathrm{rt},-1}(t)  \\
 & + b_\mathrm{hp} \dot{m}_\mathrm{hp,lt}(t) T_\mathrm{hp,lt}(t) \\
 & - \dot{m}_\mathrm{l,lt}(t) T_{\mathrm{lts},1}(t) \\
 & + \dot{m}_\mathrm{lts,+}(t) T_{\mathrm{lts},2}(t) \\
 & - \dot{m}_\mathrm{lts,-}(t) T_{\mathrm{lts},1}(t) \Bigr)
\end{aligned}\hspace{-0.4cm} \label{eq:t_lts_0_sim}
\end{equation}
where $V_\mathrm{lts}$ is a constant fluid volume of the \gls{LTS}, $T_{\mathrm{cs,w},-1}$ the \gls{CS} outlet water temperature, $T_{\mathrm{rt},-1}$ the \gls{RT} outlet temperature, $\dot{m}_\mathrm{hp,lt}$ the mass flow rate and $T_\mathrm{hp,lt}$ the outlet temperature of the \gls{LT} side of the \gls{HP} and $\dot{m}_\mathrm{l,lt}$ the \gls{LT} load mass flow. The mass flows $\dot{m}_{\mathrm{lts},+}$ and $\dot{m}_{\mathrm{lts},-}$ from and towards the neighboring storage volume, respectively, are determined by the mass balances
\begingroup
\allowdisplaybreaks
\begin{align}
\begin{split}
\dot{m}_{\mathrm{lts},+}(t) = {} & \max(0, \, \dot{m}_\mathrm{l,lt}(t) - \dot{m}_\mathrm{cs}(t) \\
& - r_\mathrm{lt,ht}(t) r_\mathrm{lts,cs}(t) \dot{m}_\mathrm{rt}(t) \\
&  - b_\mathrm{hp} \dot{m}_\mathrm{hp,lt}(t)) \label{eq:mdot_lts_plus}
\end{split} \\ 
\begin{split}
\dot{m}_{\mathrm{lts},-}(t) = {} & \max(0, \, \dot{m}_\mathrm{cs}(t) + b_\mathrm{hp} \dot{m}_\mathrm{hp,lt}(t) \\
& + r_\mathrm{lt,ht}(t) r_\mathrm{lts,cs}(t) \dot{m}_\mathrm{rt}(t) \\
& - \dot{m}_\mathrm{l,lt}(t)) \label{eq:mdot_lts_minus}
\end{split}
\end{align}
\endgroup
and are valid for all volumes within the \gls{LTS}, so it is not necessary to compute these balances for each discrete volume individually.

The discontinuities introduced to $f_\mathrm{sim}$ by \cref{eq:mdot_lts_plus,eq:mdot_lts_minus} must be avoided within $f_\mathrm{mpc}$ to preserve differentiability. One possible way to achieve this would be to use only one volume per storage, however, not representing storage stratification at all results in insufficient information for the \gls{MPC} on the current storage state and can therefore lead to severe performance losses regarding both optimality as well as constraint satisfaction.

Therefore, it should be preferred to introduce multiple volumes per storage and suitable, smooth approximations for the mass flows between those. According to the method presented by \cite{Sawant2019} and introducing a separate control variable $\dot{m}_\mathrm{rt,lts}$ for the fraction of the mass flow from the \gls{RT} through the \gls{LTS} to reduce model nonlinearity, Eq. \eqref{eq:t_lts_0_sim} is reformulated for use with\-in $f_\mathrm{mpc}$ as
\begingroup
\allowdisplaybreaks
\begin{align}
\begin{split}
\dot{m}_\mathrm{lts,s}(t) = {} & \dot{m}_\mathrm{rt,lts}(t) + \dot{m}_\mathrm{cs}(t) \\
& + b_\mathrm{hp} \dot{m}_\mathrm{hp,lt}(t) - \dot{m}_\mathrm{l,lt}(t)
\end{split} \\
\dot{m}_\mathrm{lts,\bar{s}}(t) = {} & \sqrt{\dot{m}_\mathrm{lts,s}^2(t) + \epsilon \,} \label{eq:lts_mdot_abs_approx} \\
\begin{split}
\dot{T}_{\mathrm{lts},1}(t) =  {} & (\rho_\mathrm{w} V_\mathrm{lts})^{-1} \biggl(\dot{m}_\mathrm{cs}(t) T_{\mathrm{cs,w},-1}(t)  \\
 & +  \dot{m}_\mathrm{rt,lts}(t) T_{\mathrm{rt},-1}(t) \\
 & + b_\mathrm{hp} \dot{m}_\mathrm{hp,lt}(t) T_\mathrm{hp,lt}(t) \\
 & - \dot{m}_\mathrm{l,lt}(t) T_{\mathrm{lts},1}(t) \\
 & - \Bigl( \dot{m}_\mathrm{lts,s}(t) (T_{\mathrm{lts},1}(t) + T_{\mathrm{lts},2}(t)) / 2 \\
 & + \dot{m}_\mathrm{lts,\bar{s}}(t) (T_{\mathrm{lts},1}(t) - T_{\mathrm{lts},2}(t)) / 2 \Bigr) \biggr).
\end{split}\hspace{-0.6cm} \label{eq:t_lts_0_mpc}
\end{align}
\endgroup

Though the approximation of the absolute val\-ue in \eqref{eq:lts_mdot_abs_approx} is better the smaller $\epsilon$ is chosen, smaller values for $\epsilon$ increase the nonlinearity introduced to the model, so that $\epsilon$ must be considered as a tuning parameter and chosen based on the actual quantities of the mass flows. Within this study, $\epsilon=0.5$ is used.

The \gls{HTS} is modeled analogously to the \gls{LTS}, however $T_{\mathrm{hts},1}$ marks the temperature of the top volume of the \gls{HTS}. The \gls{HR} is installed in the top and can therefore be used to directly increase $T_{\mathrm{hts},1}$. Also, the mass flow from the \gls{RT} to the \gls{HTS} is introduced as an individual control $\dot{m}_\mathrm{rt,hts}$ in $f_\mathrm{mpc}$, accordingly to $\dot{m}_\mathrm{rt,lts}$ for the \gls{LTS}.
\subsection{Heat pump}
For modeling the \gls{HP}, we assume that its internal controller realizes a constant temperature difference of $\Delta T_\mathrm{hp} = 4\,\mathrm{K}$ between inlet and outlet on both the \gls{LT} and \gls{HT} side of the \gls{HP} as in
\begin{align}
T_\mathrm{hp,ht}(t) & = T_{\mathrm{hts},n_\mathrm{hts}}(t) + \Delta T_\mathrm{hp}, \\
T_\mathrm{hp,lt}(t) & = T_{\mathrm{lts},n_\mathrm{lts}}(t) - \Delta T_\mathrm{hp}.
\end{align}

The \gls{COP} is calculated from the Carnot COP for the machine and a degree of efficiency factor $\eta_\mathrm{hp} = 0.52$ as in
\begin{equation}
\mathrm{COP}_\mathrm{hp}(t)  = \eta_\mathrm{hp} \frac{T_\mathrm{hp,ht}(t)}{T_\mathrm{hp,ht}(t) - T_\mathrm{hp,lt}(t)}.
\end{equation}

Using $\mathrm{COP}_\mathrm{hp}$ and the electrical power $P_\mathrm{hp,el}$ of the machine, the heating power $\dot{Q}_\mathrm{hp,ht}$ and cooling power $\dot{Q}_\mathrm{hp,lt}$ under current operation conditions can be computed as
\begin{align}
\dot{Q}_\mathrm{hp,ht}(t) & = \mathrm{COP}_\mathrm{hp}(t) P_\mathrm{hp,el}(t), \\
\dot{Q}_\mathrm{hp,lt}(t) & = \dot{Q}_\mathrm{hp,ht}(t) - P_\mathrm{hp,el}(t).
\end{align}

Finally, the current mass flow on the \gls{LT} and \gls{HT} side of the machine can be calculated as
\begin{align}
\dot{m}_\mathrm{hp,lt}(t) & = \frac{\dot{Q}_\mathrm{hp,lt}(t)}{c_\mathrm{w} (T_{\mathrm{lts},n_\mathrm{lts}}(t) - T_\mathrm{hp,lt}(t))}, \\
\dot{m}_\mathrm{hp,ht}(t) & = \frac{\dot{Q}_\mathrm{hp,ht}(t)}{c_\mathrm{w} (T_\mathrm{hp,ht}(t) - T_{\mathrm{hts},n_\mathrm{hts}}(t))}.
\end{align}

\subsubsection{Recooling tower}

The \gls{RT} is modeled by a number of fluid volumes that exchange heat with their environment depending on ambient temperature $T_\mathrm{amb}$ and current \gls{RT} operation level $r_\mathrm{rt,op}$. Since the \gls{RT} is designed as an air-water cross-flow heat exchanger, air mass flows are assumed to be high and resulting air temperature differences small. Following this assumption, no further air balance volume discretization is conducted. For simplicity, $r_\mathrm{tr,op}$ is modeled to directly influence the current heat exchange of the \gls{RT}, so that the temperature $T_\mathrm{rt}$ within a balance volume is calculated as
\begin{align}
\dot{T}_\mathrm{rt}(t) = {} & \frac{\dot{m}_\mathrm{rt}(t) (T_\mathrm{rt,pre}(t) - T_\mathrm{rt}(t)) - \frac{\dot{Q}_\mathrm{rt}(t)}{c_\mathrm{w}}}{\rho_\mathrm{w} V_\mathrm{rt}} \\
\dot{Q}_\mathrm{rt}(t) = {} & r_\mathrm{rt,op}(t) \alpha_\mathrm{rt} A_\mathrm{rt} (T_\mathrm{rt}(t) - T_\mathrm{amb}(t))
\end{align}
with $V_\mathrm{rt}$ the constant fluid volume and $A_\mathrm{rt}$ the heat exchange area between fluid and outside environment depending on the discretization. For the volume at the inlet of the \gls{RT}, $T_\mathrm{rt,pre}$ results from the mixing of the outlet temperatures of \gls{CS}, \gls{LTS} and \gls{HTS} according to their mass flow ratios, for the upcoming volumes it is the temperature of the preceding volume.

While for the \gls{MPC} model the \gls{RT} is modeled using a single volume, the \gls{RT} is discretized into $n_\mathrm{rt} = 5$ volumes for the simulation model.

\subsection{Modeling of thermal loads of the building}

Cooling loads $\dot{Q}_\mathrm{l,lt}$ are covered from the \gls{LTS} and heating loads $\dot{Q}_\mathrm{l,ht}$ from the \gls{HTS}. In the following, we assume that an internal, low-level controller realizes a constant temperature difference $\Delta T_\mathrm{l}$ by supply flow addition in the heating and cooling circuits, while setting the mass flows in the distribution circuits to meet demands. The return temperatures and mass flows entering and leaving the storages are then given by
\begin{align}
T_\mathrm{l,lt}(t) = {} & T_{\mathrm{lts},1}(t) + \Delta T_\mathrm{l}, \\
\dot{m}_\mathrm{l,lt}(t) = {} &  \left(c_\mathrm{w} \Delta T_\mathrm{l} \right)^{-1} \dot{Q}_\mathrm{l,lt}(t), \\
T_\mathrm{l,ht}(t) = {} & T_\mathrm{hts,1}(t) - \Delta T_\mathrm{l}, \\
\dot{m}_\mathrm{l,ht}(t) = {} & \left(c_\mathrm{w} \Delta T_\mathrm{l} \right)^{-1} \dot{Q}_\mathrm{l,ht}(t) .
\label{eq:loads}
\end{align}

While utilization of a detailed building model may facilitate further efficiency gains within the later control applications, this would also further increase the complexity of the model through addition of multiple states and control variables, and is therefore  considered future work.

\subsection{Auxiliary states}

The \gls{MPC} model contains the auxiliary states $C_{\dot{Q}_\mathrm{l,lt}}$, $C_{\dot{Q}_\mathrm{l,ht}}$ and $C_{T_\mathrm{amb}}$ used for compensation of deviations between forecasted and measured time-varying parameters to increase performance of the \gls{MPC}. Those are determined by solving an \gls{MHE} problem in between two subsequent \gls{MPC} steps. Exemplary, the auxiliary state for $T_\mathrm{amb}$ reads as
\begin{equation}
\dot{C}_{T_\mathrm{amb}}(t) = -(\tau_{C}^{-1} C_{T_\mathrm{amb}}(t))
\end{equation}
with $\tau_C$ a time constant which allows for modeling a decaying influence of the current deviation with increasing time. With this and $T_\mathrm{amb,fc}$ the raw forecast inputs, $T_\mathrm{amb}$ then enters $f_\mathrm{mpc}$ as
\begin{equation}
T_\mathrm{amb}(t) = T_\mathrm{amb,fc}(t) + C_{T_\mathrm{amb}}(t).
\end{equation}

\section{Development of a conventional control strategy}
\label{sec:conventional_control}

\begin{figure}[tbp]
\centering
\includegraphics[scale=0.35]{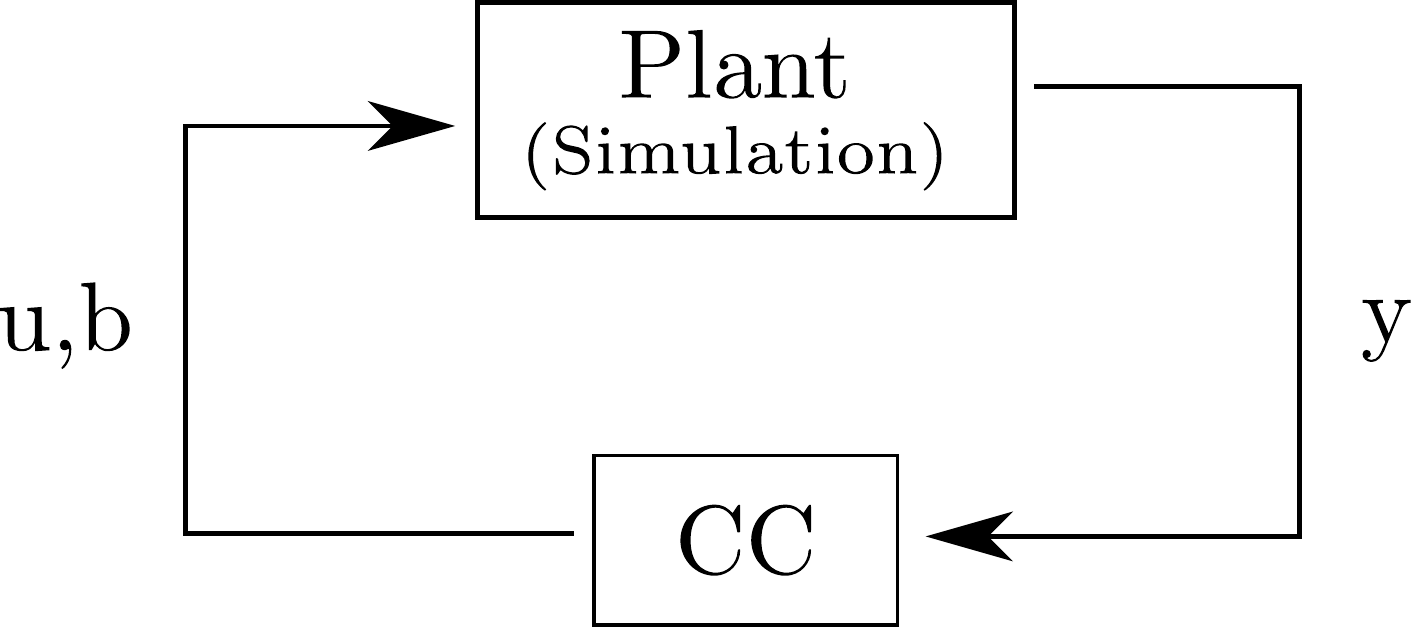}
\caption{Schematic depiction of the CC loop.}
\label{fig:cc_loop}
\end{figure}

In the following, an elaborate \gls{CC} strategy for the system is described, which has been developed and improved iteratively using the presented simulation model. A schematic depiction of the \gls{CC} loop is given in Figure~\ref{fig:cc_loop}.

\subsection{General controller design principles}

The \gls{CC} strategy aims at covering heating and cooling demands in a reliable and simultaneously energy efficient manner, realized through a sophisticated storage management approach. Since requirements regarding storage management differ between seasons, separate operation modes have been developed for summer and winter operation, which are described in more detail in the upcoming sections. Transition between summer and winter mode is carried out at fixed time points during the year.

Activation of devices is mainly realized using set-point-based control and switching hystereses depending on current temperatures of \gls{LTS}, \gls{HTS}, \gls{RT} and \gls{CS}. Mass flows $\dot{m}_\mathrm{rt}$, $\dot{m}_\mathrm{cs}$ and ratios $r_\mathrm{lts,cs}$, $r_\mathrm{rt,op}$ are determined using PI-controllers depending on specified temperature differences between inlets and outlet of the several components. Despite technically possible, $r_\mathrm{lt,ht}$ is not chosen on a continuous scale within the \gls{CC} but always fully switched depending on whether the \gls{HP} is currently mainly used for covering of cooling or heating demands.

\subsection{Operation in summer mode}

Summer mode is especially focused on covering cooling loads with irregularly occurring heating demands. For this, the following hierarchy for activation of cooling devices for decreasing \gls{LTS} temperature is defined: \begin{enumerate*}[label=\arabic*)]
\item{free cooling via the \gls{RT},}
\item{storing energy from the cooling loads in the \gls{CS}, and}
\item{cooling via the \gls{HP}.}
\end{enumerate*}

In times of spare cooling capacity, the \gls{LTS} and the \gls{CS} are cooled down using free cooling via the \gls{RT}. In case their temperatures have exceeded a certain threshold value, also the \gls{HP} is used to pro-actively cool down these storages in order to generate cooling energy under more favorable conditions and store it for later usage. If the \gls{HTS} temperature increases too much during \gls{HP} operation, the \gls{HTS} is cooled down via the \gls{RT}. Occurring heating loads are covered utilizing the \gls{HP}.

In case the \gls{HP} is active, the current driving power $P_\mathrm{hp,el}$ for the \gls{HP} is modulated in order to cover the maximum of \begin{enumerate*}[label=\alph*)]
\item{the current heating load and}
\item{the current cooling load that cannot be covered using \gls{RT} and \gls{CS}.}
\end{enumerate*}

\subsection{Operation in winter mode}

Winter mode is especially focused on covering heating loads while accounting for the permanently occurring cooling demands caused by the running server computers, but also for irregularly occurring cooling demands arising, e.\,g., in case of high solar irradiation on the building. Free cooling via the \gls{RT} is avoided and storing the energy from cooling loads in the \gls{CS} is preferred to facilitate later usage for covering of heating demands via the \gls{HP}.

In case the heating demands exceed the cooling demands in winter mode, \gls{RT} and \gls{CS} are utilized to increase the \gls{LTS} temperature to prevent undercooling. In case the \gls{LT} side of the system cannot provide enough energy for operation of the \gls{HP}, the \glspl{HR} are utilized to cover remaining loads.

\subsection{Setup and implementation}

Simulation model $f_\mathrm{sim}$ and the \gls{CC} strategy are implemented in Modelica and simulated using Dymola 2019 with DASSL integrator, which directly facilitates a simulation of the \gls{CC} strategy.

\section{Development of a model-predictive control strategy}
\label{sec:mpc_strategy}

\begin{figure}[tbp]
\centering
\includegraphics[scale=0.35]{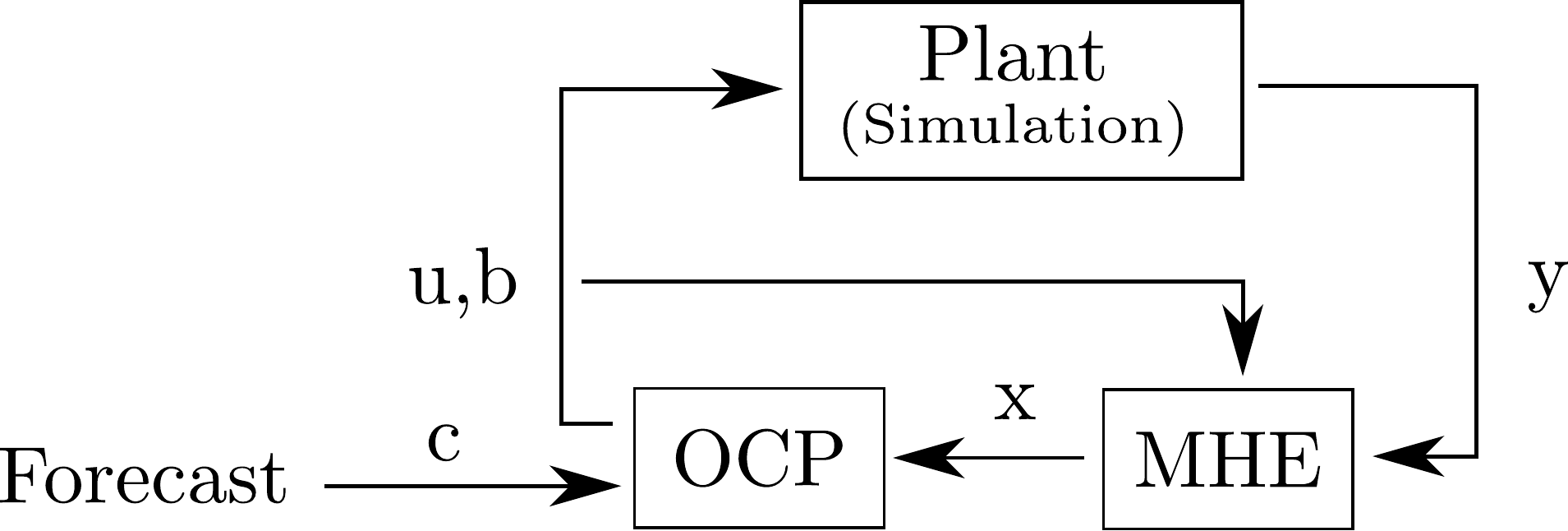}
\caption{Schematic depiction of the MPC loop.}
\label{fig:mpc_loop}
\end{figure}

In this section, the \gls{MPC} strategy is presented by first introducing the \gls{OCP} formulated for the controller and the utilized solution strategy and software implementation, followed by the formulation and implementation of the \gls{MHE} problem used for compensation of model mismatch and forecast deviations. A schematic depiction of the \gls{MPC} loop is given in Figure~\ref{fig:mpc_loop}.

\subsection{OCP formulation}

Within the \gls{OCP}, all boundaries on non-auxi\-li\-ary system states, which are the temperature constraints for the several components, are formulated as soft constraints to prevent infeasibility of the \gls{OCP} in case of state constraint violations, cf. \cite{Rawlings2017}. For this, we introduce a vector of $n_s = 33$ slack variables $s$.

With continuous controls $u^{\top} = \lbrack \dot{m}_\mathrm{rt,hts} ~ \allowbreak  \dot{m}_\mathrm{rt,lts} ~  \allowbreak  \dot{m}_\mathrm{rt,cs} ~ \allowbreak \dot{m}_\mathrm{cs} ~ \allowbreak r_\mathrm{rt,op} ~ \allowbreak P_\mathrm{hr} ~ \allowbreak  P_\mathrm{hp,el} \rbrack$ and time-varying parameters $c^{\top} = \lbrack T_\mathrm{amb,fc} ~ \allowbreak \dot{Q}_\mathrm{l,lt,fc} ~ \allowbreak \dot{Q}_\mathrm{l,ht,fc} \rbrack$, the \gls{MIOCP} reads as 
\begingroup
\allowdisplaybreaks
\begin{subequations}
\label{eq:hvac_ocp}
\begin{alignat}{2}
&\backup\begin{aligned}
\underset{{\substack{x(\cdot),u(\cdot),\\ b_\mathrm{hp}(\cdot),s(\cdot)}}}{\mathrm{min}} & \hspace{0.2cm} \int_{t_{0}}^{t_\text{f}} \left( s(t)^\top W_s s(t) + w_s^\top s(t) \right) \, \text{d}t  \\
& \hspace{0.2cm} +  \int_{t_{0}}^{t_\text{f}} w_u^\top u(t) \, \text{d}t \\
\end{aligned}\label{eq:hvac_ocp_objective} \\
\mathrm{s.\,t.} & \hspace{0.1cm} \text{for } t \in [t_0, t_\mathrm{f}]: \notag\\
& \hspace{0.1cm} \dot{x}(t) = f_\mathrm{mpc}(x(t), u(t), b_\mathrm{hp}(t), c(t)),\label{eq:hvac_ocp_dynamics} \\
& \hspace{0.1cm} 0 \leq \sum \dot{m}_\mathrm{rt,\{hts, lts, cs\}}(t) \leq \dot{m}_\mathrm{rt,max}, \label{eq:hvac_ocp_mdot_max_rt} \\
& \hspace{0.1cm} 0 \leq \dot{m}_\mathrm{rt,cs}(t) + \dot{m}_\mathrm{cs}(t) \leq \dot{m}_\mathrm{cs,max}, \label{eq:hvac_ocp_mdot_max_cs} \\
& \hspace{0.1cm} T_\mathrm{lb} - s(t) \leq Z x(t) \leq T_\mathrm{ub} + s(t), \label{eq:hvac_ocp_states_limits_bounds} \\
& \hspace{0.1cm} u_\mathrm{lb} \leq u(t) \leq u_\mathrm{ub}, \label{eq:hvac_ocp_controls_limits_bounds} \\
& \hspace{0.1cm} b_\mathrm{hp}(t) \in \{0, 1\}, \label{eq:hvac_ocp_limits_binary} \\
& \hspace{0.1cm} s(t) \geq 0, \label{eq:hvac_ocp_limits_slacks} \\
& \hspace{0.1cm} x(t_0) = x_{0}. \label{eq:hvac_ocp_initial_state}
\end{alignat}
\end{subequations}
\endgroup

The objective of the \gls{MIOCP} given in \eqref{eq:hvac_ocp_objective} is a Lagrange-type objective and contains the sum of squares of the slack variables, the sum of the slack variables and the sum of the continuous controls, weighted by appropriate diagonal weighting matrix $W_\mathrm{s} \in \mathbb{R}^{n_s \times n_s}$ and vectors $w_\mathrm{s} \in \mathbb{R}^{n_s}$ and $w_\mathrm{u} \in \mathbb{R}^{n_u}$, respectively. The dynamics of the system are given in \eqref{eq:hvac_ocp_dynamics}. Equations \cref{eq:hvac_ocp_mdot_max_rt,eq:hvac_ocp_mdot_max_cs} ensure that the maximum mass flows through \gls{RT} and \gls{CS} are not exceeded. With $Z$ a selection matrix that selects the constrained temperatures from $x$, the soft constraints for the state boundaries are given in \eqref{eq:hvac_ocp_states_limits_bounds}, while the boundaries on the continuous controls are given as hard constraints in \eqref{eq:hvac_ocp_controls_limits_bounds}. Finally, the binary constraint is given in \eqref{eq:hvac_ocp_limits_binary}, positivity of the slack variables is ensured by \eqref{eq:hvac_ocp_limits_slacks} and the initial state constraint is given in \eqref{eq:hvac_ocp_initial_state}.

\subsection{OCP setup and implementation}
\label{subsec:ocp_setup_and_implementation}

For solving \glspl{OCP}, utilization of direct methods, especially direct multiple shooting \citep{Bock1984} or direct collocation \citep{Tsang1975}, is favorable, cf. \cite{Sager2009,Binder2001}. For implementation and solution of \eqref{eq:hvac_ocp}, a control horizon of $t_\mathrm{f} = 24\,\mathrm{h}$ has been chosen and discretized into $N=29$ control intervals, whereof the first 6 intervals are of length $\Delta t_\mathrm{s} = 10\,\mathrm{min}$ and the remaining 23 intervals of length $\Delta t_\mathrm{l} = 1\,\mathrm{h}$, using direct collocation with Lagrange polynomials with Radau collocation points.

For solution of the resulting \gls{MINLP} on \gls{MPC} suitable time scales, application of general \gls{MINLP} solvers is usually not favorable. Therefore, the decomposition approach by \cite{Sager2009} is applied, where the solution of the original \gls{MINLP} is approximated by solving first a relaxed version of the \gls{MINLP}, which is a \gls{NLP}, then obtaining an integer trajectory for the discrete controls by application of either \gls{SUR} \citep{Sager2009} or \gls{CIA} \citep{Sager2011a}, and afterwards solving the \gls{NLP} again with binary controls fixed according to the solution of the integer approximation step. In the literature, several successful applications of this approach can be found \citep{kirches2011fast,Ebrahim2018,Buerger2018}.

The \gls{NLP} is implemented using the dynamic optimization framework CasADi \citep{Andersson2018} in Python. For solution of the arising \glspl{NLP}, Ipopt \citep{Waechter2006} with linear sol\-ver MA57 \citep{HSL2017} is used. For solution of the integer approximation problem, the binary approximation package pycombina is applied, which is presented in \cite{Buerger2019}.

The controls for the system obtained by solution of the \gls{MINLP} for the first 10 min of the control horizon are applied to the simulation model $f_\mathrm{sim}$ implemented in Modelica to accept external control inputs via the Python interface of Dymola.

\begin{figure*}[tbp]
\centering
\includegraphics[width=\textwidth]{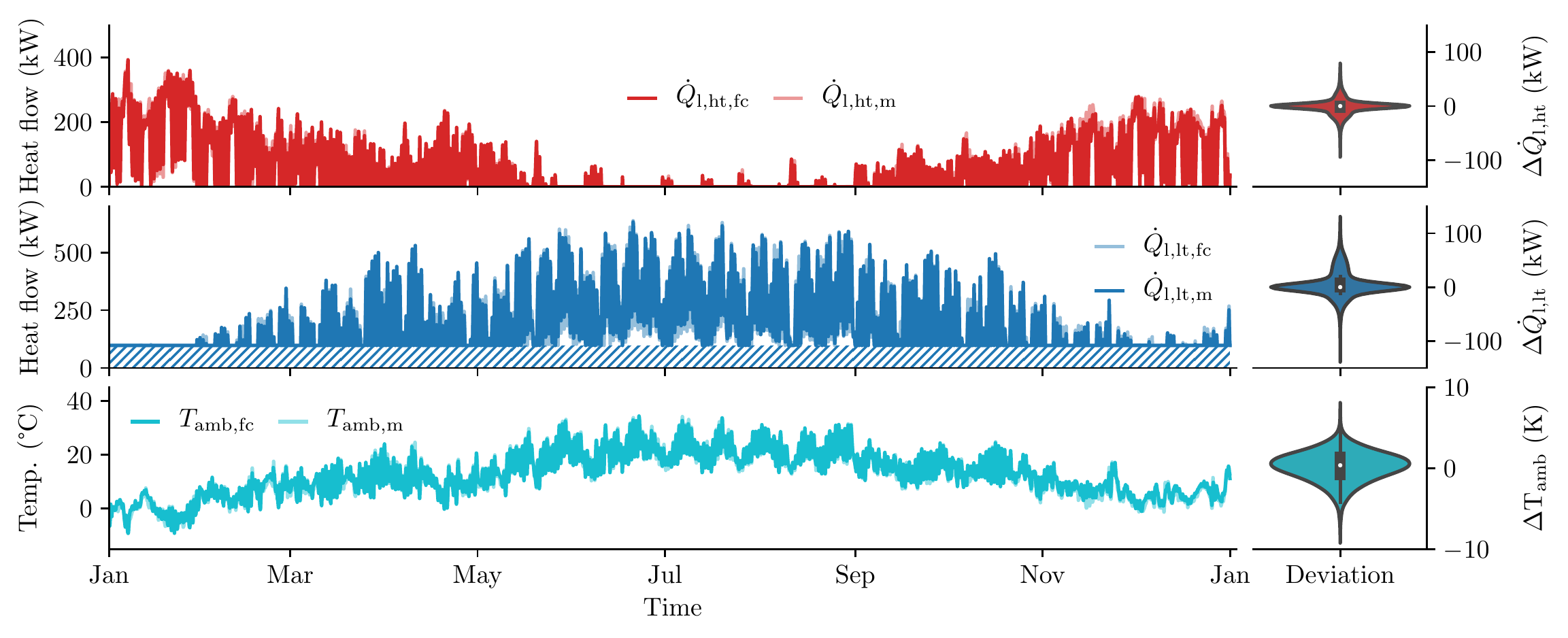}
\caption{Heating demands $\dot{Q}_\mathrm{l,lt}$, cooling demands $\dot{Q}_\mathrm{l,ht}$ and ambient temperatures $T_\mathrm{amb}$ assumed within this study. In the left hand side plots, suffix ``fc'' indicates forecasts used for the \gls{MPC} and suffix ``m'' indicates corresponding measured values used for simulation. The hatched blue area marks the constant cooling demand induced by the server computers. The right hand side plots are violin plots showing the distribution of the deviations from forecasts to measurements.}
\label{fig:params}
\end{figure*}

\subsection{MHE formulation}

The \gls{MHE} problem for estimation of the current system state and compensation of the current forecast deviations based on the previous $M=4$ measurements reads as
\begingroup
\allowdisplaybreaks
\begin{subequations}
\label{eq:hvac_mhe}
\begin{alignat}{2}
&\backup\begin{aligned}
\underset{x(\cdot),w(\cdot)}{\mathrm{min}} & \hspace{0.2cm} \sum_{k = 0}^{M} \| y(x_k, u_k, b_{\mathrm{hp},k}, c_k) - x_k \|_{W^{-1}}^2 \\
& \hspace{0.2cm} + \sum_{k = 0}^{M} \| w_k \|_{Q^{-1}}^2 + \| x_0 - \hat{x} \|_{P^{-1}}^2 \\
\end{aligned}\label{eq:hvac_mhe_objective} \\
\mathrm{s.\,t.} & \hspace{0.1cm} \text{for } k = 0, \dots, M-1: \notag\\
&\hspace{0.2cm} x_{k+1} = \bar{f}_\mathrm{mpc}(x_k, u_k, b_{\mathrm{hp},k}, c_k, w_k) \label{eq:hvac_mhe_dynamics}
\end{alignat}
\end{subequations}
\endgroup
with $y$ the measurement function and $W$ the covariance matrix of the measurements, $w$ an additive process noise with covariance matrix $Q$ acting on the system state, $\| x_0 - \hat{x} \|_{P^{-1}}^2$ the arrival cost with covariance matrix $P$, and $\bar{f}_\mathrm{mpc}$ the solution of the \gls{MPC} model for a discrete time step.

More detailed information on \gls{MHE} and interpretation of the arrival cost can be found in, e.\,g., \cite{Kuehl2011} and \cite{Kraus2013}.

\subsection{MHE setup and implementation}

Similar to the \gls{OCP}, the \gls{MHE} problem is implemented using direct collocation and the arising \glspl{NLP} are solved using Ipopt and MA57. An application of the decomposition approach is not necessary, as the \gls{MHE} problem contains no discrete optimization variables. Arrival cost updates are conducted using a smooth\-ened \gls{EKF}, cf. \cite{Girrbach2018}.

The resulting measurements after simulating the first 10 min of the horizon in Dymola are passed as measurements to the \gls{MHE}. The estimated system states resulting from solution of the \gls{MHE} problem are then used for initialization of the subsequent \gls{MPC} step. By repeating this procedure, a simulation of the \gls{MPC} strategy is carried out.

\section{Comparison of controller performances}
\label{sec:controller_performance_comparison}

In the following, the controller performances are compared on an annual simulation for each controller. All computations are conducted on a Fujitsu P920 Desktop PC with an Intel Core i5-4570 3.20 GHz CPU and 16 GB RAM running Debian 9, using Python 3.5, gcc 6.3, Dymola 2019, CasADi 3.4.5 and pycombina 0.2.

\subsection{Scenario of the simulation study}
The forecasted and measured thermal loads $\dot{Q}_{\mathrm{l,ht},\{\mathrm{fc,m}\}}$ and $\dot{Q}_{\mathrm{l,lt},\{\mathrm{fc,m}\}}$ of the building assumed within this study are calculated using a dynamic building simulation in accordance to VDI 2078 based on forecasts and measured weather data provided by \cite{KAweather2017}. These are depicted in Figure~\ref{fig:params} next to the forecasted and measured ambient temperature $T_{\mathrm{amb},\{\mathrm{fc,m}\}}$. 

For realistic initialization of the system states, their initial values have been chosen based on the final states obtained from preliminary simulations.

\subsection{Comparison of simulation results}

\begin{figure*}[tbp]
\centering
\includegraphics[width=\textwidth]{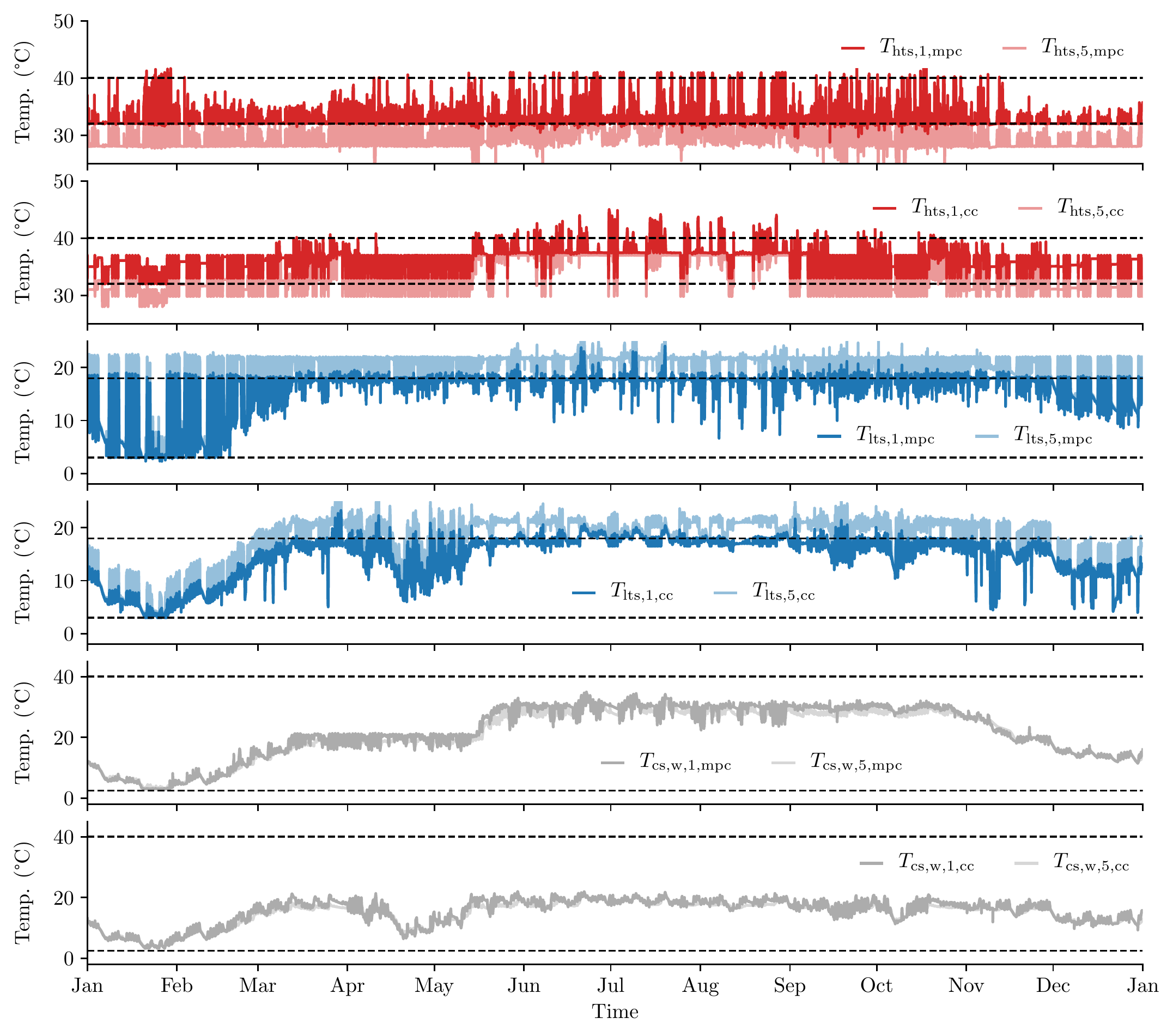}
\caption{Simulated temperatures within \gls{HTS}, \gls{LTS} and \gls{CS} throughout the year obtained by \gls{MPC} and \gls{CC}. The upper plots show the temperatures $T_\mathrm{hts,1}$ at the top and $T_\mathrm{hts,5}$ at the bottom of the \gls{HTS}, the dashed black lines indicate the temperature constraints for $T_\mathrm{hts,1}$. The middle plots show the temperatures $T_\mathrm{lts,1}$ at the bottom and $T_\mathrm{lts,5}$ at the top of the \gls{LTS}, the dashed black lines indicate the temperature constraints for $T_\mathrm{lts,1}$. The lower plots show the water temperatures $T_\mathrm{cs,w,1}$ at the inlet and $T_\mathrm{cs,w,5}$ at the outlet of the \gls{CS} and the dashed black lines denote their constraints.}
\label{fig:storages}
\end{figure*}

\begin{figure*}[t]
\centering
\includegraphics[width=\textwidth]{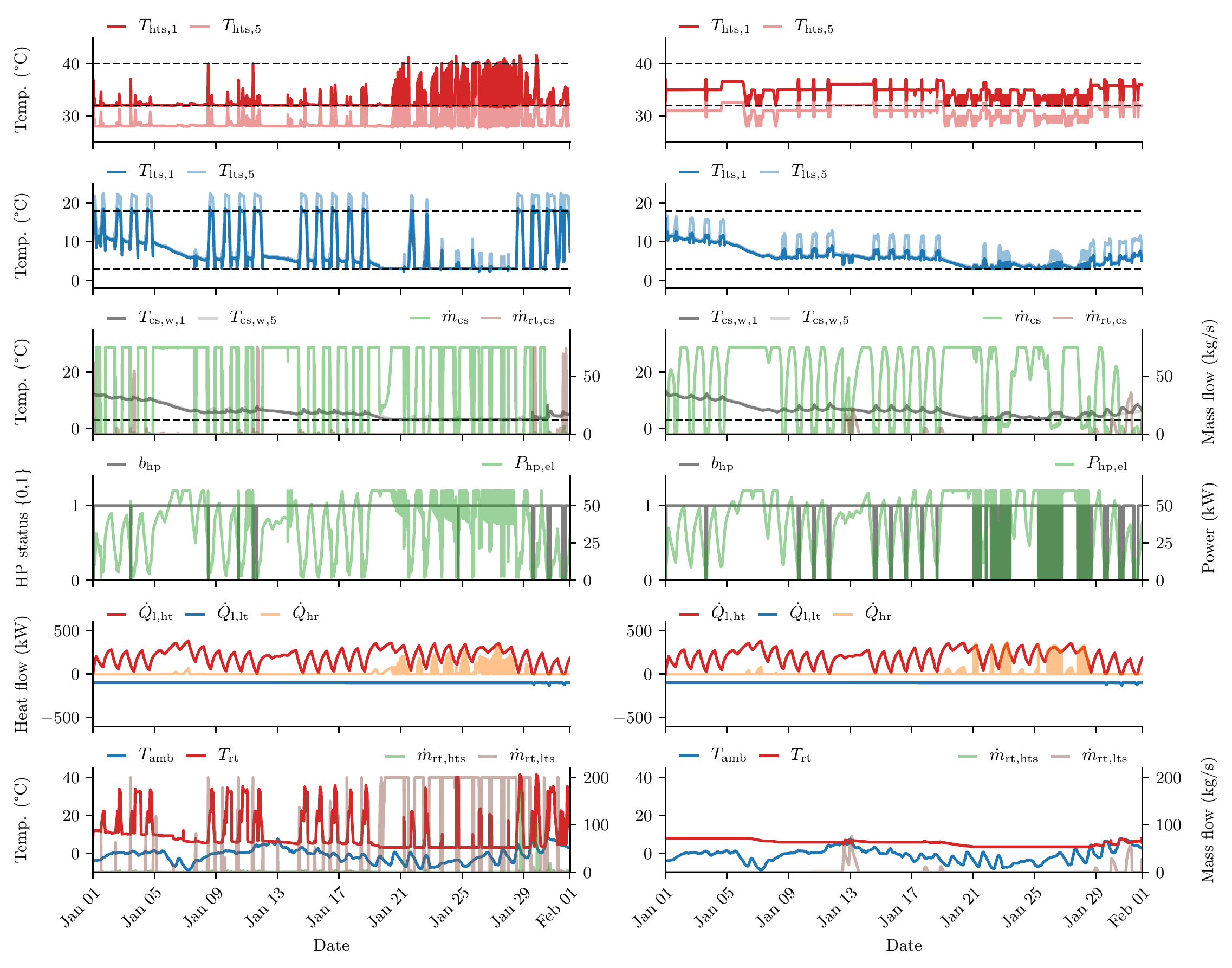}
\caption{Detailed view on system operation using \gls{MPC} (left) and \gls{CC} (right) in January.}
\label{fig:results_detail_jan}
\end{figure*}

An overview of the results of the simulation study is given in Figure~\ref{fig:storages} by the temperatures of the \gls{HTS}, \gls{LTS} and \gls{CS} throughout the year for both \gls{MPC} and \gls{CC}. The plot reveals substantially different behaviors of the controllers, which is discussed in more detail in the following.

\subsubsection{LTS and HTS management}

Figure~\ref{fig:storages} shows that the \gls{CC} aims to achieve a certain margin of the storage temperatures from their boundary values. This is necessary for the controller to stay reactive in case of altering loads. Therefore, it takes action once a temperature val\-ue moves close towards a boundary and drives the temperature value towards the opposite direction.

For the \gls{MPC} however, it can be observed that the controller runs the storage on temperatures which are very close to their respected boundaries. Explicitly, the \gls{MPC} tries to always run the \gls{LTS} at the upper temperature boundary and the \gls{HTS} at the lower temperature boundary. This increases the efficiency of the \gls{HP}, and consequently, decreases its electrical power consumption.

\begin{figure*}[t]
\centering
\includegraphics[width=\textwidth]{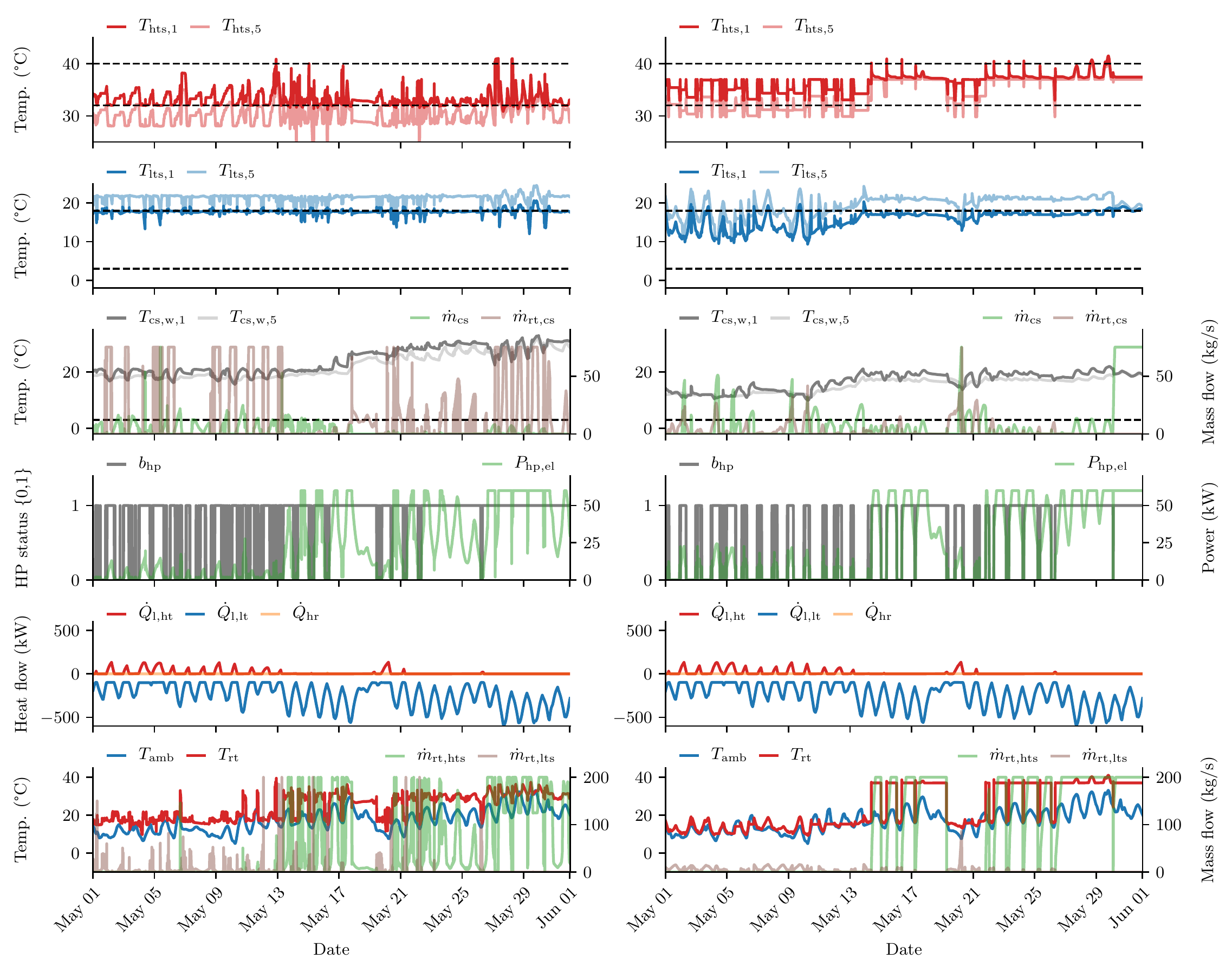}
\caption{Detailed view on system operation using \gls{MPC} (left) and \gls{CC} (right) in May.}
\label{fig:results_detail_may}
\end{figure*}

A more detailed depiction of this behavior can be found in Figure~\ref{fig:results_detail_jan}, which depicts the trajectories of states and controls in January. Here, the \gls{LTS} acts as a heat source for the \gls{HP}, and it can be observed how the \gls{MPC} favors \gls{HP} efficiency by increasing the temperature of the \gls{LTS} periodically in times of low heating load while keeping the \gls{HTS} temperature at its lower boundary. Consistently, it can be seen that the electric power $P_\mathrm{hp,el}$ consumed by the \gls{HP} is often lower for the \gls{MPC} than for the \gls{CC}.

Operation of the storages this close to their respected boundary temperatures is only possible due to the incorporation of (corrected) forecasts in the \gls{MPC} strategy and its knowledge of the system behavior through the included model $f_\mathrm{mpc}$. This is a clear advantage of the \gls{MPC} over the \gls{CC}, however, coming at the price of an increased complexity of controller design and implementation.

\subsubsection{Concrete slab utilization}

\begin{figure*}[t]
\centering
\includegraphics[width=\textwidth]{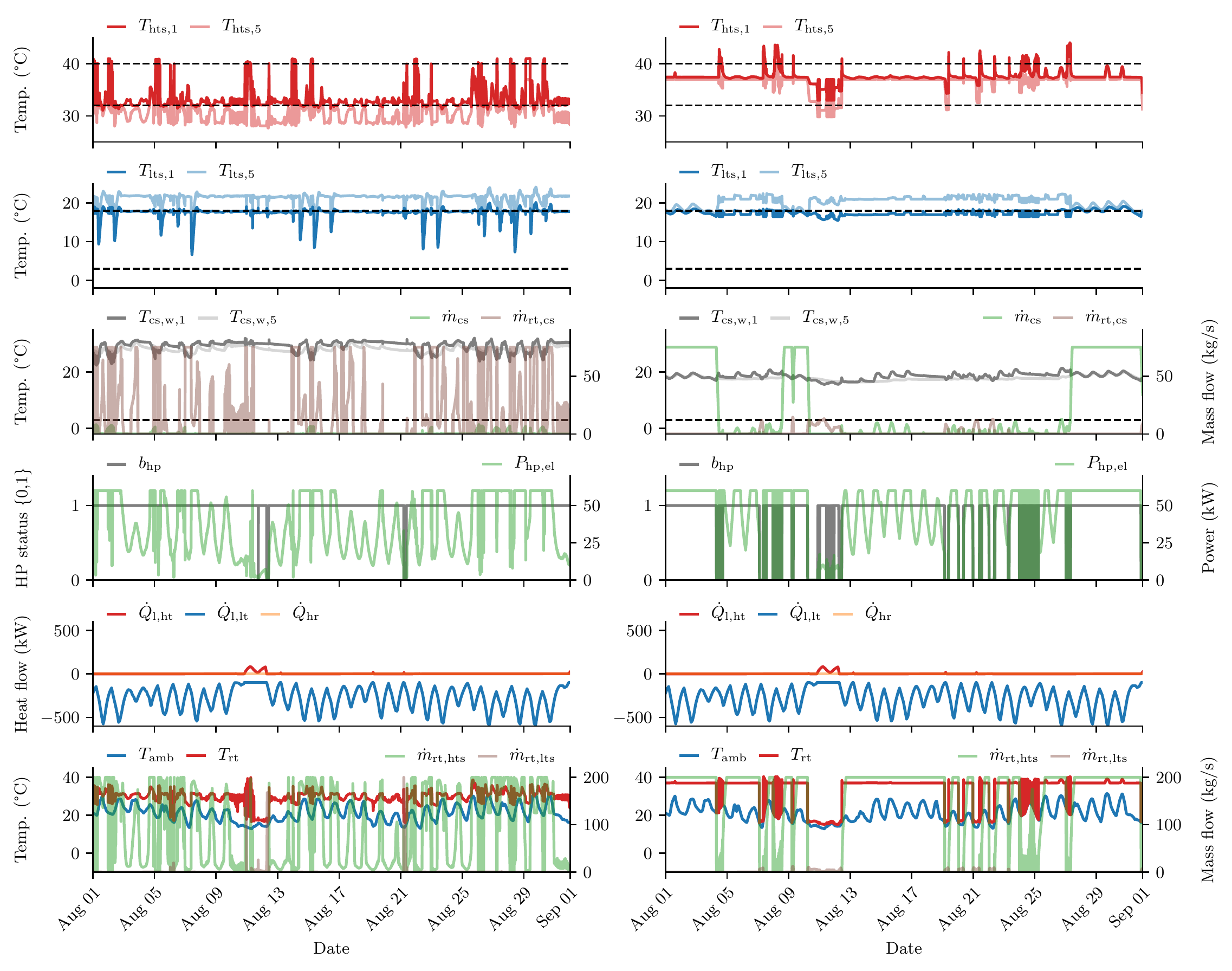}
\caption{Detailed view on system operation using \gls{MPC} (left) and \gls{CC} (right) in August.}
\label{fig:results_detail_aug}
\end{figure*}

The winter mode of the \gls{CC} strategy, which is designed for times dominated by heating loads, utilizes the \gls{CS} to increase the temperature of the \gls{LTS}, and with this, as a heat source for the \gls{HP} covering the occurring heating loads. During such load scenarios, the \gls{MPC} utilizes the \gls{CS} similarly to the \gls{CC}, as shown by the \gls{CS} temperatures in Figure~\ref{fig:storages} and by the \gls{CS} temperatures, pump operation and \gls{HP} operation in Figure~\ref{fig:results_detail_jan}.

In times dominated by cooling loads though, the strategies of the controllers differ. Here, the summer mode of the \gls{CC} utilizes the \gls{CS} to provide cooling power for reducing the temperature of the \gls{LTS}, i.\,e., to store energy from the \gls{LTS}, and is regenerated in times of low ambient temperatures using the \gls{RT} or, in rare cases, the \gls{HP}.

The \gls{MPC} however utilizes the \gls{CS} to store energy also from the \gls{HTS}. More insight into the reason behind this can be found in Figure~\ref{fig:results_detail_may}, which depicts the trajectories of states and controls in May. While it is shown that in the first half of the month the \gls{MPC} uses the \gls{CS} similarly to the \gls{CC} as a cold source for the \gls{LTS}, in the second half of the month the \gls{MPC} realizes mass flows from the \gls{HTS} to the \gls{RT} and from the \gls{CS} to the \gls{RT} at the same time, which results in the \gls{CS} being used as a recooling device for the \gls{HP}. Due to the lower recooling temperature supported to the \gls{HTS} in comparison of using only the \gls{RT} individually, the operation efficiency of the \gls{HP} increases.

The energy stored in the \gls{CS} this way can then be emitted to the ambient via the \gls{RT} during times of lower ambient temperature. A higher temperature difference between \gls{CS} and ambient can also be achieve easier now as if the \gls{CS} was used to merely cool down the \gls{LTS} which operates at a much lower maximum temperature than the \gls{HTS}. Further, this results in a higher utilization of the \gls{CS} within the \gls{MPC} when compared to the \gls{CC}, which can be observed in Figure~\ref{fig:results_detail_aug}, which depicts the trajectories of states and controls in August.

This utilization of the \gls{CS} for recooling of the \gls{HP}, which turns out to be an efficient and favorable operation mode, was identified autonomously by the \gls{MPC} but not during the design phase of the \gls{CC}, which shows the advantages of \gls{MPC} for flexible identification of situation-dependent and non-intuitive operation modes.

\subsubsection{Constraint satisfaction}

\begin{figure*}[t]
\centering
\includegraphics[width=\textwidth]{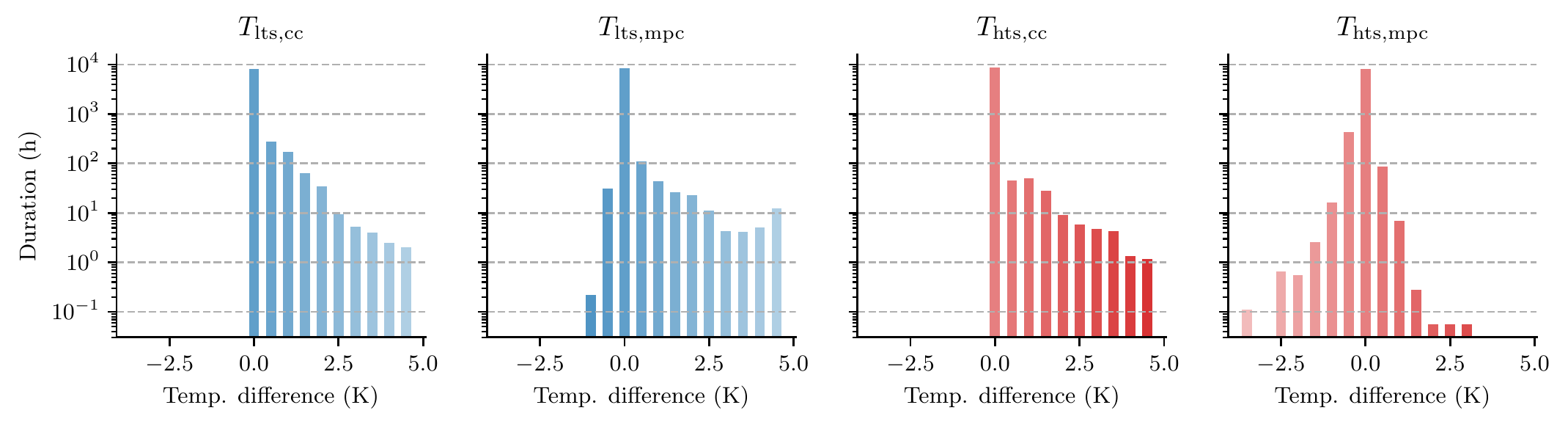}
\caption{Total duration and severity of storage temperature constraint violations over the whole simulation time.}
\label{fig:constraint_violation}
\end{figure*}

Figure~\ref{fig:storages} already indicates that for the relevant storage temperatures $T_\mathrm{hts,1}$ and $T_\mathrm{lts,1}$ as well as for the temperatures of the \gls{CS}, both controllers achieve a high quality of constraint satisfaction. A more detailed analysis on the duration and levels of constraint violations throughout the year is shown in Figure~\ref{fig:constraint_violation}.

The plot shows the total duration and severity of the differences between these storage temperatures and their upper and lower boundaries over the whole simulation time. For both storage temperatures, their respective boundaries are violated by more than 0.5\,K in less than 572\,h throughout the year, which corresponds to less than 7\,\% of the simulated time.

For the \gls{CC}, one can observe that the lower temperature boundaries are never violated. This is due to the implementation of the \gls{CC}, which does not operate on a fixed time grid and directly stops \gls{HP} operation in the event of too low \gls{LTS} temperatures and activates the \glspl{HR} in the event of too low \gls{HTS} temperatures.

The \gls{MPC}, however, operates on a fixed time grid that allows changes of the requested operations of components only every $\Delta t_\mathrm{s} = 10\,\mathrm{min}$. Since the \gls{MPC} tends to operate the storage close to their boundary values, sudden mismatches between forecasted and occurring loads can cause violations of these boundaries. Such violations could possibly be further reduced by decreasing the \gls{MPC} time step size $\Delta t_\mathrm{s}$, and consequently, the reaction time on occurring disturbances. However, due to the already high computation time for deduction of this \gls{MPC} simulation study (cf. Section~\ref{subsec:runtime}), this has not been further investigated.

While the violations of the upper temperature boundaries for the \gls{LTS} are comparable for \gls{CC} and \gls{MPC}, the upper boundaries for the \gls{HTS} tend to be violated less for the \gls{MPC} than for the \gls{CC}, which can be explained by the \gls{MPC} generally operating the \gls{HTS} at its lower boundary.

\begin{figure*}[tbp]
\centering
\includegraphics[width=\textwidth]{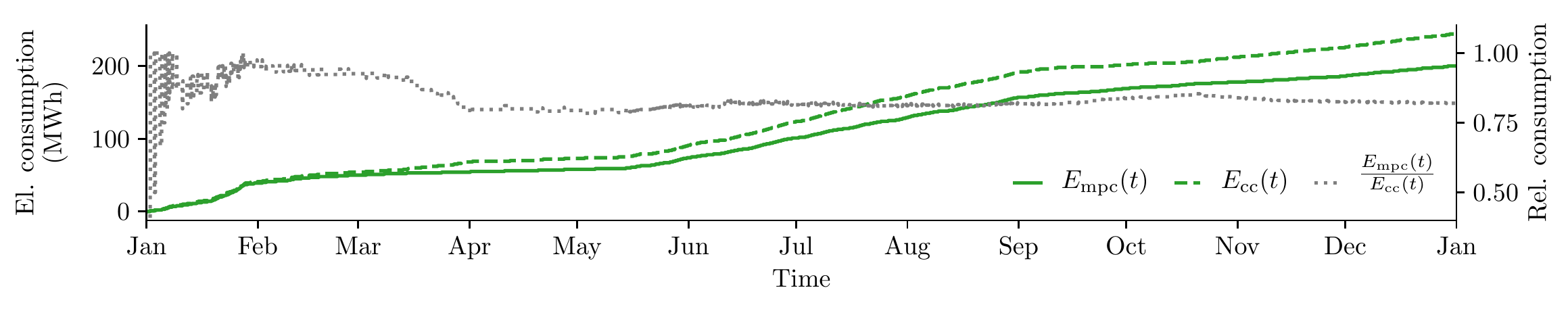}
\caption{Integrated electrical energy consumption $E(t) = \int_{t_0}^{t} \bigl( P_\mathrm{hp,el}(\tau) + P_\mathrm{hr}(\tau) \bigr) \mathrm{d}\tau$ of \gls{HP} and \glspl{HR} for \gls{MPC} and \gls{CC}.}
\label{fig:consumption}
\end{figure*}

\subsubsection{Energy consumption}

The integrated electrical energy consumption of the \gls{HP} and the \glspl{HR}\footnotemark are shown for each controller in Figure~\ref{fig:consumption}. Until April, the consumptions are comparable and only slightly lower for the \gls{MPC}, which can be explained by the similar behavior of both controllers until that time and the relatively small number of options for control actions to satisfy the high heating loads during these times.

\footnotetext{Within this work, only the electrical energy consumption of the \gls{HP} and \glspl{HR}, which are regarded the main electrical consumers of the system, are considered, since the exact pressure drops and corresponding pump dimensions as well as the exact \gls{RT} configuration, which determine the electrical consumption of these components, were not defined for the system at the time of this study. Considerations on the electricity consumption of pumps could provide further potential for operation optimization, cf. \cite{Wystrcil2013}.}

After that, the relative consumption for the \gls{MPC} reduces throughout the rest of the year due to the more efficient operation modes identified by the \gls{MPC} during summer times. Also, the resulting higher amount of energy stored in the \gls{CS} at the end of the warmer season can be used for driving the \gls{HP} more efficiently in upcoming colder seasons. At the end of the year, the main electrical energy consumption is reduced by more than 18\,\% using \gls{MPC}.

Additionally, using \gls{MPC} drastically reduces the amount of power-cycles for the \gls{HP}, which is switched merely 2135 times over the year, compared to 7075 power-cycles requested by the \gls{CC}.

\begin{figure}[tbp]
\centering
\includegraphics[width=\columnwidth]{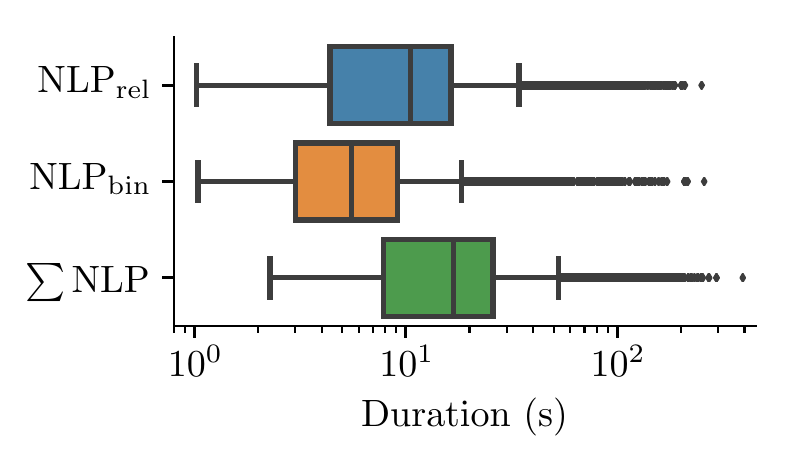}
\caption{Durations of the runtime-critical stages of the \gls{MINLP} solution process for all N = 52,560 steps. The boxes range from the 25\,\% to the 75\,\% quartile, the whiskers of the boxes extend these by 1.5 times the inter-quartile range. Values outside of this range are regarded as outliers and denoted as dots.}
\label{fig:durations}
\end{figure}

\subsection{Runtime of critical \gls{MINLP} solution steps}
\label{subsec:runtime}

Figure~\ref{fig:durations} shows the runtimes of those \gls{MINLP} solution steps with significant contribution to the overall runtime of the algorithm, which are the solutions of the \glspl{NLP} described in Section~\ref{subsec:ocp_setup_and_implementation}. The runtime of the binary approximation step and the solution of the \gls{MHE} problem within this study is always within the range of seconds and therefore negligible.

The plot shows that both the typical solution time as well as the maximum duration of both steps are below 400\,s and with this far below the step size of the \gls{MPC}, which renders the method generally applicable also for real-time usage on the physical system.

\subsection{Significance of the \gls{MHE} step}
\label{subsec:mhe_impact}

To investigate significance and impact of the \gls{MHE} step on the \gls{MPC} performance, several attempts have been made to run the simulation study without a dedicated state estimator and the auxiliary states described above, only utilizing simple approaches for correction of forecasts, e.\,g., calculation of a rolling mean on the deviations between measured and forecasted inputs. However, none of these approaches was able to sufficiently correct the forecasts, which at some point during simulation caused the \gls{MPC} to fail due to severe storage temperature constraint violations. Prior to failing, constraint satisfaction was generally worse and small to medium sized violations occurred more often compared to using \gls{MHE}. These findings render the \gls{MHE} step significant for the overall \gls{MPC} performance.

\section{Conclusion}
\label{sec:conclusion}

We presented a whole-year simulation study on nonlinear mixed-integer \gls{MPC} for a complex thermal energy supply system and compared its performance to a conventional, set point based control strategy. We showed that \gls{MPC} was able to reduce the yearly electrical energy consumption of \gls{HP} and \glspl{HR} by more than 18\,\% while providing a similar degree of constraint satisfaction. Furthermore, the \gls{MPC} enabled the planner to identify previously unknown, beneficial operation modes.

Future work should consider incorporation of a detailed building model, which could facilitate further efficiency gains. Considerations on the electricity consumption of pumps could provide further potential for operation optimization. Also, through consideration of dynamic electricity pri\-ces, the system could be utilized in the context of demand response.

\section*{Conflict of interest}

None.

\section*{Acknowledgments}

The authors kindly thank Ottensmeier Ingenieure Gmbh, Paderborn, Germany for their cooperation and support of this project.

This research was supported by the German Federal Ministry for Economic Affairs and Energy (BMWi) via eco4wind (0324125B) and DyConPV (0324166B), by DFG via Research Unit FOR 2401, by the State Ministry of Baden-Wuert\-temberg for Sciences, Research and Arts (Az: 22-7533.-30-20/9/3) and by the German Federal Ministry for the Environment, Nature Conservation, Building and Nuclear Safety (BMUB) via WIN4\-climate (0204KF0354).

\printglossary[style=tree, type=\acronymtype, title=Acronyms, nonumberlist]

\bibliographystyle{elsarticle-harv} 
\bibliography{literature}

\end{document}